\documentclass{elsarticle}

\usepackage{hyperref}
\usepackage{float}
\usepackage{graphicx}
\usepackage{subfigure}
\usepackage{amsmath,bm}
\usepackage{epstopdf} 
\usepackage{caption}
\usepackage{natbib}
\biboptions{sort&compress}

\usepackage{mathptmx}
\usepackage{verbatim}
\usepackage{indentfirst}
\usepackage[title]{appendix}
\usepackage{amssymb}
\usepackage{cases}
\usepackage{color}

\begin{document}
\newtheorem{theorem}{Theorem}[section]
\newtheorem{example}{Example}[section]
\newtheorem{remark}{Remark}[section]
\numberwithin{equation}{section}
\numberwithin{figure}{section}
\numberwithin{table}{section}

\captionsetup[figure]{font=small,labelfont={bf},labelformat={default},labelsep=period,name={Fig.}}
\captionsetup[table]{labelfont={bf},labelformat={default},labelsep=period,name={Table}}

\begin{frontmatter}

\title{Macroscopic auxiliary asymptotic preserving neural networks for the linear radiative transfer equations}

\author[10614]{Hongyan Li}
\ead{lihongyan142869@163.com}

\author[909]{Song Jiang}
\ead{jiang@iapcm.ac.cn}

\author[909]{Wenjun Sun}
\ead{sun$\_$wenjun@iapcm.ac.cn}

\author[10614]{Liwei Xu\corref{mycorrespondingauthor}}
\cortext[mycorrespondingauthor]{Corresponding author}
\ead{xul@uestc.edu.cn}

\author[106141]{Guanyu Zhou}
\ead{zhoug@uestc.edu.cn}

\address[10614]{School of Mathematical Sciences, University of Electronic Science and Technology of China, Chengdu Sichuan 611731, China}
\address[909]{Institute of Applied Physics and Computational Mathematics, Beijing 100094, China}
\address[106141]{Institute of Fundamental and Frontier Sciences, University of Electronic Science and Technology of China, Chengdu Sichuan 611731, China}

\begin{abstract}
We develop a Macroscopic Auxiliary Asymptotic-Preserving Neural Network (MA-APNN) method to solve the time-dependent linear radiative transfer equations (LRTEs), which have a multi-scale nature and high dimensionality.
To achieve this, we utilize the Physics-Informed Neural Networks (PINNs) framework and design a new adaptive exponentially weighted Asymptotic-Preserving (AP) loss function, which incorporates the macroscopic auxiliary equation that is derived from the original transfer equation directly and explicitly contains the information of the diffusion limit equation.
Thus, as the scale parameter tends to zero, the loss function gradually transitions from the transport state to the diffusion limit state. 
In addition, the initial data, boundary conditions, and conservation laws serve as the regularization terms for the loss. We present several numerical examples to demonstrate the effectiveness of MA-APNNs.
\end{abstract}


\end{frontmatter}


\section{Introduction}
The kinetic equations describe the motion of particles in a medium, including collision and absorption. These equations are defined in phase space and can involve multiple spatial and temporal scales, as well as nonlocal operators, making numerical simulations difficult\cite{bouchut2000kinetic,dimarco2014numerical,jin2010asymptotic,jin2023asymptotic1}. LRTEs, as a typical kinetic model, find wide-ranging applications in fields such as astrophysics, weapon physics, inertial/magnetic confinement fusion\cite{chandrasekhar2013radiative,larsen1983asymptotic,adams1998asymptotic}, and more. Accurately simulating LRTEs is nontrivial, mainly due to two challenges. The first challenge is high dimensionality. LRTEs describe the evolution of the radiation density function of a vast number of photons, which, in the most general cases, is a seven-dimensional problem. When uncertainties are factored into the equations, the dimensionality increases further\cite{jin2017uncertainty,hu2016stochastic,jin2017uniform}. The second challenge is the multiscale features caused by different opacities of the background materials. LRTEs typically involve multiple spatial and/or temporal scales, characterized by the Knudsen number, which is the dimensionless mean free path. When the Knudsen number is tiny, LRTEs can be approximated by the diffusion limit equations\cite{bouchut2000kinetic}.
To preserve the propagation properties of photons in different optical regions, the multiscale modeling approaches are introduced to the numerical simulation\cite{weinan2011principles,degond2005smooth}.
When coupled radiation transport and diffusion models are solved using the same spatial grid, the mesh size should be roughly equivalent to the mean free path. This requirement can make computations incredibly expensive in optically thick mediums where the Knudsen number is very small.

Efficient computational methods that can deal with high dimensionality and multi-scale characteristics are highly desirable. Classical tools such as the Implicit Monte Carlo (IMC) method \cite{fleck1971implicit, gentile2001implicit, mcclarren2009modified, shi2020continuous, densmore2011asymptotic} are capable of handling high dimensionality, but the convergence rate is low, and there are statistical noises. On the other hand, Asymptotic-Preserving (AP) schemes, which are uniformly stable with respect to the scaling parameter, are proposed to address the multiscale problem. As the limit of the Knudsen number approaches zero, the solution of the AP schemes converges to the diffusion limit case.
The AP schemes were initially proposed to solve the neutron transport problems by Larsen and Keller\cite{larsen1974asymptotic}, Larsen\cite{larsen1974solutions,larsen1975neutron,larsen1976neutron}, and then improved by Larsen,
Morel and Miller \cite{larsen1987asymptotic}, Larsen and Morel\cite{morel1989asymptotic}, Jin and Levermore\cite{jin1991discrete,jin1993fully}.
For the unsteady cases, \cite{jin2000uniformly,klar1998asymptotic} also constructed the AP schemes by decomposing the distribution into the equilibrium and disturbance parts.
These methods were further developed and extended to other multiscale kinetic models\cite{jin1999efficient,jin2010asymptotic,mieussens2013asymptotic,lemou2008new}.

In recent years, deep learning methods and deep neural networks (DNNs) have been developed vigorously and achieved some success in solving PDEs \cite{han2018solving,zhang2021mod,weinan2021algorithms,hwang2020trend,beck2020overview}.
Various approaches have been proposed, including the Physics-Informed Neural Networks (PINNs) \cite{raissi2019physics} and the deep Galerkin method (DGM) \cite{sirignano2018dgm} based on the least-squares formulation,
the deep Ritz method (DRM) \cite{yu2018deep} utilizing the variational formulation, and others \cite{zang2020weak,cai2021least,lyu2022mim}. Compared to traditional mesh-based methods, the DNNs have an advantage as mesh-free methods in solving problems with complex geometric domains and high dimensions \cite{poggio2017and,grohs2018proof}. Furthermore, DNNs use automatic differentiation, which avoids truncation errors in discretization. However, DNN methods have some drawbacks, including lack of high accuracy, long training time due to a large number of network parameters, and the use of stochastic gradient method for optimizing high-dimensional non-convex functions\cite{lu2021deepxde}.

The PINNs, Model-Operator-Data Network (MOD-Net), and Physics-Informed DeepONets (PIDONs) have been applied to solve steady and unsteady linear radiative transfer models\cite{zhang2021mod,mishra2021physics,lu2022solving,wang2021learning,chen2021solving,liu2021deep}. For MOD-Net method, the DNN is used to parameterize the Greens function and the neural operator is obtained to approximate the solution. Combining the operator universal approximation theory and PINN's idea, the PIDON method can solve a class of PDEs from the learned continuous operators. We find that the vanilla PINN, PIDON, or MOD-Net approach, has difficulties dealing with the multiscale characteristics. The loss function of these methods deteriorates as the scale parameter decreases to zero, which fails to capture asymptotic limits and usually tends to learn a simplified model during the training process. To tackle this problem, the Asymptotic-Preserving Neural Networks (APNNs)\cite{jin2023asymptotic1} and Model-Data APNNs\cite{liSongSunXuZhou2023} are designed by using the micro-macro decomposition to solve the time-dependent linear transport equations and gray radiative transfer equations, respectively.
In addition, \cite{jin2023asymptotic2} develops the APNN method based on the even-odd decomposition for the multiscale kinetic equations. The micro-macro and even-odd decomposition techniques have also been applied to develop the Asymptotic-Preserving Convolutional Deep Operator Networks (APCONs) for the time-dependent linear radiative transport problem\cite{wu2023asymptotic}.

In this work, we devise a novel loss function incorporating the macroscopic auxiliary equation, which is capable of capturing the diffusion limit behavior as the scale parameter approaches zero. We derive the macroscopic auxiliary equation for the LRTEs by repeatedly substituting the original equation into itself and integrating it with respect to the angle direction. This was inspired by the idea of deriving the diffusion limit system for the radiative heat transport equation as described in \cite{ghattassi2022diffusive}. The macroscopic auxiliary equation we derived holds true for smooth solutions and can be considered a perturbation of the diffusion limit equation when the scale parameter is small.
In addition, we introduce an exponential weight that depends on the scale parameter. This weight is used to combine the original equation's loss with that of the macroscopic auxiliary equation, to achieve the AP property. We also consider the mass conservation constraints of the radiation intensity under periodic boundaries, as well as the boundary and initial hard constraints under inflow and periodic boundaries. To evaluate the effectiveness of the proposed Macroscopic Auxiliary APNN (MA-APNN) method, we simulate low and high-dimensional linear transport equations under transport and diffusion regimes. We investigate various scenarios, such as initial layer problems, constant or variable scattering problems, uncertainty quantification problems, and more. Compared to PINNs\cite{raissi2019physics}, MA-APNNs have great advantages in solving diffusion dominating problems. Different from the APNNs in \cite{jin2023asymptotic1,jin2023asymptotic2} where the macro-micro or even-odd decomposition is used, it does not need the decomposition of radiation intensity in MA-APNNs, so the network's construction is simple, and in turn cause the training time is reduced greatly with the comparable prediction accuracy as APNNs\cite{jin2023asymptotic1} under diffusion case.

The rest of this paper is organized as follows. In Section \ref{sec2}, we introduce the time-dependent linear radiative transfer model and its corresponding diffusion limit equation. In Section \ref{sec3}, we present the motivation of the MA-APNN method and define the new AP loss function. Numerous numerical experiments, including both multiscale and high dimensional uncertainty problems, are carried out in Section \ref{sec4} to validate the efficiency of MA-APNNs. A conclusion remark is given in Section \ref{sec5}.

\section{Linear radiative transfer equations}\label{sec2}
The linear radiative transfer equations describe the radiation photons transport and the energy exchange with the background materials.
We denote by $f(t,r,\Omega)$ the radiation intensity of photons located at the space point $r=(x,y,z) \in D$ in time $t \in \tau$ and propagating in direction in $\Omega = (\xi,\eta,\mu) \in S^2$, where $D$ is a bounded domain in $\mathbb{R}^3$ and $S^2$ is the unit sphere. The scaled form of the LRTEs with initial and boundary conditions is stated as follows.
\begin{subequations}\label{eq:LRTEs}
	\begin{align}
		& \epsilon^2 \frac{\partial f}{\partial t}+\epsilon \Omega \cdot \nabla_r f =\sigma \left(\frac{1}{4 \pi}\int_{S^2}f\text{d}\Omega-f\right) - \epsilon^2 \alpha f + \epsilon^2G && \text{ in } \tau\times D \times S^2,  \label{eq:LRTEs-a} \\
		& B f = f_b && \text{ on } \tau\times \partial D \times S^2,  \label{eq:LRTEs-b}\\
		& f(0,r,\Omega) = f_0(r,\Omega) && \text{ in } D \times S^2, \label{eq:LRTEs-c}
	\end{align}
\end{subequations}
where $\sigma(r)$ is the scattering coefficient, $\alpha(r)$ the absorption coefficient, $G(r)$ the internal source, and $\epsilon >0$ the Knudsen number.
The average of $f(t,r,\Omega)$ on $S^2$ is called the incident radiation $\rho(t,r)$, saying
\[
\rho := \left\langle f \right\rangle = \frac{1}{4 \pi}\int_{S^2}f\text{d}\Omega.
\]
Integrating \eqref{eq:LRTEs-a} with respect to $\Omega$, we obtain
\begin{equation}\label{eq:LRTEs-a-rho}
\epsilon^2 \frac{\partial \rho}{\partial t}+\epsilon \langle \Omega \cdot \nabla_r f \rangle = - \epsilon^2 \alpha \rho + \epsilon^2  G  \quad \text{ in } \tau\times D.
\end{equation}
Furthermore, integrating \eqref{eq:LRTEs-a-rho} with respect to $r$ and applying the integration by parts, we get the mass conservation equation:
\begin{equation}\label{eq:LRTEs-mass}
\epsilon^2 \frac{d}{dt} \int_D \rho\text{d}r + \epsilon \frac{1}{4\pi} \int_{\partial D} \int_{S^2} \Omega \cdot n_r f\text{d}\Omega \text{d}r = - \epsilon^2 \alpha \int_D \rho\text{d}r + \epsilon^2  \int_{D} G\text{d}r,
\end{equation}
where $n_r$ represents the unit outer normal to $\partial D$ at position $r$.
In the following, we state the boundary condition \eqref{eq:LRTEs-b} in detail, and introduce the 1D and 2D models, the uncertainty quantification model, and the diffusion limit equation.
\subsection{The boundary conditions}
We separate the boundary $\Gamma= \tau\times\partial D \times S^2$ into the inflow and outflow parts as follows:
\[
\Gamma_- := \left\{(t,r,\Omega)\in \Gamma : \Omega \cdot n_{r} <0\right\}, \qquad \Gamma_+ := \left\{(t,r,\Omega)\in \Gamma : \Omega \cdot n_{r} >0\right\}.
\]
In this work, we consider two types of boundary conditions for LRTEs.
\begin{enumerate}[(1)]
	\item The inflow boundary condition
	\begin{equation}\label{eq:BC-Dir}
	f(t,r,\Omega) = f_{B}(t,r,\Omega) \quad (t,r,\Omega)\in \Gamma_-,
	\end{equation}
	where $f_{B}(t,r,\Omega)$ is a given function.
	\item If $D$ is symmetric, we can enforce the periodic boundary condition \cite{lemou2008new}
	\begin{equation}\label{eq:BC-per}
	f(t,r,\Omega)=f(t, \mathcal{S}(r),\Omega) \quad (t,r,\Omega) \in \Gamma,
	\end{equation}
	where $\mathcal{S}:\partial D \rightarrow \partial D$ is a one-to-one symmetry mapping satisfying $n_r = - n_{\mathcal{S}(r)}$.
	At this time, we see from \eqref{eq:LRTEs-mass} that
	\begin{equation}\label{eq:mass-conservation}
		\frac{d}{dt}\int_{D}\rho \text{d}r + \int_{D}\alpha\rho \text{d}r - \int_{D} G \text{d}r = 0.
	\end{equation}
\end{enumerate}
\subsection{The 1D and 2D linear models}
In one-dimensional case, i.e., $(r,\Omega)=(x, \mu) \in (x_L, x_R) \times (-1,1)$, the LRTEs becomes
\begin{equation}\label{eq:LRTEs-1d-lin}
\epsilon^2 \partial_t f+\epsilon \mu \partial_x f =\sigma \left(\frac{1}{2}\int_{-1}^{1}f \text{d} \mu -f\right) - \epsilon^2\alpha f  + \epsilon^2 G,
\end{equation}
with the isotropic inflow boundary condition
\begin{equation}\label{eq:LRTEs-1d-inflow-bc}
	f(t,x_L,\mu>0)=f_L(t,\mu), \quad f(t,x_R,\mu<0)=f_R(t,\mu),
\end{equation}
or the periodic boundary condition
\begin{equation}\label{eq:LRTEs-1d-periodic-bc}
	f(t,x_L,\mu)=f(t,x_R,\mu),
\end{equation}
and the initial condition
\begin{equation}\label{eq:LRTEs-1d-initial}
	f(0,x,\mu)=f_0(x,\mu).
\end{equation}
In two-dimensional case, i.e., $(r,\Omega) = ((x,y),(\mu,\xi)) \in D \times S^1$, \eqref{eq:LRTEs-a} is replaced by
\begin{equation}\label{eq:LRTEs-2d-lin}
\epsilon^2 \partial_t f+\epsilon \Omega \cdot \nabla_r f  =\sigma \left(\frac{1}{2\pi}\int_{S^1}f \text{d} {\bm{v}} -f\right) - \epsilon^2\alpha f  + \epsilon^2 G.
\end{equation}
where $D\subset \mathbb{R}^2$ and $S^1=\{(\xi,\eta) : \xi^2+\eta^2=1\}$.
\subsection{ The uncertainty quantification (UQ) model}
In practice, we are interested in the LRTEs that contains uncertainty in the collision cross section, source, initial or boundary data\cite{jin2017uniform}.
The uncertainty is characterized by the random variable
$\bm{z} =(\bm{z}_1, \bm{z}_2,...,\bm{z}_m) \in \mathbb{R}^m$ with the probability density $\omega(\bm{z})$.
In our numerical simulation, we consider the spatial 1D LRTEs with the random input $\bm{z}$, where the radiation intensity $f(t,x,\mu,\bm{z})$ satisfies the equation
\begin{equation}\label{eq:LRTEs-uq}
\epsilon^2 \partial_t f+\epsilon \mu \partial_x f =\sigma(\bm{z}) \left(\frac{1}{2}\int_{-1}^{1}f \text{d} \mu -f\right) - \epsilon^2\alpha(\bm{z}) f  + \epsilon^2 G(\bm{z}),
\end{equation}
and the initial and inflow boundary conditions
\begin{subequations}\label{eq:LRTEs-uq-initial}
\begin{align}
& f(0,x,\mu, \bm{z})=f_0(x,\mu, \bm{z}), \label{eq:LRTEs-uq-initial-a} \\
& f(t,x_L,\mu, \bm{z})=f_L(t,\mu,\bm{z}) \quad  (\mu>0), \label{eq:LRTEs-uq-initial-b}\\
& f(t,x_R,\mu, \bm{z})=f_R(t,\mu, \bm{z}) \quad  (\mu<0). \label{eq:LRTEs-uq-initial-c}
\end{align}
\end{subequations}

\subsection{The diffusion limit equation}
As the Knudsen number $\epsilon \rightarrow 0$, the radiation intensity $f(t,r,\Omega)$ converges to the average density $\rho=\left\langle f \right\rangle$ and satisfies the asymptotic diffusion limit equation\cite{jin2000uniformly}, i.e.,
\begin{subequations}\label{eq:Diff-lim-all}
\begin{align}
& f = \rho , \label{eq:Diff-lim-f} \\
& \partial_t \rho - \langle \Omega^2 \rangle \nabla_r\left(\frac{1}{\sigma}\nabla_r\rho\right) + \alpha\rho - G =0, \label{eq:Diff-lim}
\end{align}
\end{subequations}
where
\[
\langle \Omega^2 \rangle := \frac{1}{4\pi}\int_{S^2} \mu^2 \text{d}\Omega = \frac{1}{3}\frac{1}{4\pi}\int_{S^2} \mu^2 +\xi^2+\eta^2  \text{d}\Omega = \frac{1}{3}\frac{1}{4\pi}\int_{S^2} 1  \text{d}\Omega = \frac{1}{3}.
\]
In 1D case, $\langle \Omega^2 \rangle := \frac{1}{2} \int_{-1}^1 \mu^2\text{d}\mu =1/3$, while in 2D case, $\langle \Omega^2\rangle:=\frac{1}{2\pi} \int_{S^1} \xi^2 \text{d}\Omega = 1/2$. An asymptotic-preserving scheme for \eqref{eq:LRTEs} should be uniformly stable with respect to $\epsilon$ and lead to an accurate approximation to the diffusion limit equation \eqref{eq:Diff-lim} when $\epsilon$ approaches zero.

\section {Macroscopic auxiliary asymptotic preserving neural networks}\label{sec3}
We first introduce the notations for DNNs. Given an input $\bm{x}=(t,r,\Omega) \in \tau \times D\times S^2$, the output of the $L$-layer feedforward neural network $f_{\theta}$ is defined recursively as
\begin{equation}
\begin{aligned}{\label{eq:DNN}}
&f_{\theta}^{[0]}(\bm{x})=\bm{x}, \\
&f_{\theta}^{[l]}(\bm{x})=\sigma^{h} \circ \left(W^{[l-1]}f_{\theta}^{[l-1]}(\bm{x})+b^{[l-1]}\right), \quad l\le l\le L-1, \\
&f_{\theta}(\bm{x})=f_{\theta}^{[L]}(\bm{x})=\sigma^{o} \circ \left(W^{[L-1]}f_{\theta}^{[L-1]}(\bm{x})+b^{[L-1]}\right),
\end{aligned}
\end{equation}
where $W^{[l]}\in \mathbb{R}^{m_{l+1}\times m_l}$ and $b^{[l]}\in \mathbb{R}^{m_{l+1}}$ represent weight matrix and bias vector, respectively, $\sigma^{h}$ is the nonlinear activation function of the hidden layers, and $\sigma^{o}$ is a specially designed activation function of the output layer. The common choices of $\sigma^{h}$ are the sigmoid function, the hyperbolic tangent function, the ReLU function, and so on. The notation ``$\circ$'' means entry-wise composition. For simplicity, we denote the set of parameters by $\theta$ and represent the network by a list, i.e., $[m_0, m_1,...,m_L]$, where $m_0$ and $m_L$ are the dimensions of the input and the output.

Solving PDE using DNNs involves three critical steps. Firstly, we parameterize the solution of PDE with a deep neural network. Secondly, we design the population/empirical loss function associated with the PDE, which evaluates the error between the approximate and exact solutions. Finally, we select appropriate optimization algorithms to minimize the loss function.

According to \cite{jin2023asymptotic1}, for the multiscale problem, it is crucial to build a loss function that preserves the asymptotic property. In the following, we first briefly explain why the vanilla PINNs fail to resolve the multiscale LRTEs. Then, we present a formal derivation of the macroscopic auxiliary equation. After that, we build a loss function incorporating the macroscopic auxiliary equation and demonstrate the asymptotic-preserving property.

\subsection{The vanilla PINNs fail to resolve the multiscale LRTEs}\label{eq:subsec-PINNs}
Let us consider the LRTEs with periodic boundary condition \eqref{eq:BC-per}, which satisfies the mass conservation law \eqref{eq:mass-conservation}.
To parametrize the nonnegative radiation intensity $f(t,r,\Omega)$, we use the neural network $f^{nn}_{\theta}(t,r,\Omega)$ with the nonnegative activation function $\sigma^o(\bm{x})=e^{-\bm{x}}$ at the output layer,
and adopt the least square of the residual of LRTEs as the loss function:
\[
L_{\text{PINNs}}^{\epsilon} = L_{\text{PINNs},g}^{\epsilon} + L_{\text{PINNs},b}^{\epsilon} + L_{\text{PINNs},i}^{\epsilon} + L_{\text{PINNs},c}^{\epsilon},
\]
where
\[
\begin{aligned}
& L_{\text{PINNs},g}^{\epsilon} = \lambda_g \left\| \epsilon^2\partial_t f_{\theta}^{nn} + \epsilon \Omega\cdot \nabla_r f_{\theta}^{nn}- \sigma\left(\left \langle f_{\theta}^{nn} \right \rangle - f_{\theta}^{nn} \right) + \epsilon^2\alpha f_{\theta}^{nn} - \epsilon^2 G \right\|^2_{L^2(\tau \times D \times S^2)}, \\
& L_{\text{PINNs},b}^{\epsilon}  = \lambda_{b} \left\| B f_{\theta}^{nn} - f_b \right\|^2_{L^2(\tau \times \partial D \times S^2)}, \\
& L_{\text{PINNs},i}^{\epsilon}  = \lambda_{i}  \left\| f_{\theta}^{nn}(0) - f_0 \right\|^2_{L^2(D \times S^2)}, \\
& L_{\text{PINNs},c}^{\epsilon} = \lambda_{c} \Big \Vert \epsilon \partial_t\int_{D}\left \langle f_{\theta}^{nn} \right \rangle \text{d}r + \int_{\partial D} \langle \Omega \cdot n_r f_{\theta}^{nn} \rangle \text{d}r  + \epsilon \int_{D}\alpha \left \langle f_{\theta}^{nn} \right \rangle \text{d}r - \epsilon \int_{D} G \text{d}r  \Big \Vert^2_{L^2(\tau)}.
\end{aligned}
\]
Here, $\lambda_g, \lambda_i, \lambda_b$ and $\lambda_c$ are the weight parameter to be tuned. Note that $L_{\text{PINNs},c}^{\epsilon}$ is the residual of the mass conservation constraint.
If the periodic boundary condition is imposed, then, by \eqref{eq:mass-conservation},
\[
L_{\text{PINNs},c}^{\epsilon} = \lambda_{c} \Big \Vert \partial_t\int_{D}\left \langle f_{\theta}^{nn} \right \rangle \text{d}r + \int_{D}\alpha \left \langle f_{\theta}^{nn} \right \rangle \text{d}r - \int_{D} G \text{d}r  \Big \Vert^2_{L^2(\tau)}.
\]
When $\epsilon$ is not too small, one can expect to obtain the approximate solution by minimizing the above loss function $L_{\text{PINNs}}^{\epsilon}$.
However, as $\epsilon \rightarrow 0$, we find that
\[
L_{\text{PINNs},g}^{\epsilon}\rightarrow \lambda_g \left\| - \sigma
		\left(\left \langle f_{\theta}^{nn} \right \rangle - f_{\theta}^{nn} \right) \right\|^2_{L^2(\tau \times D \times S^2)},
\]
which is nothing but the residual of the following equation
\[
\sigma (\langle f \rangle - f) = \sigma (\rho - f) = 0.
\]
Apparently, it is not the desired diffusion limit equation \eqref{eq:Diff-lim}.

If we add the residual of the diffusion limit equation
\[
L_{\text{PINNs},d}^{\epsilon} = \lambda_d \left\| \partial_t \left \langle f_{\theta}^{nn} \right \rangle - \left\langle \Omega^2 \right \rangle \nabla_r\left(\frac{1}{\sigma}\nabla_r \left \langle f_{\theta}^{nn} \right \rangle \right) + \alpha \left \langle f_{\theta}^{nn} \right \rangle - G \right\|^2_{L^2(\tau \times D)},
\]
to the total loss directly, saying
\[
L_{\text{PINNs}}^{\epsilon} = L_{\text{PINNs},g}^{\epsilon} + L_{\text{PINNs},b}^{\epsilon} + L_{\text{PINNs},i}^{\epsilon} + L_{\text{PINNs},c}^{\epsilon} + L_{\text{PINNs},d}^{\epsilon},
\]
then we will recover the diffusion limit system as $\epsilon \rightarrow 0$.
But, unfortunately, it fails to solve the problem with large or medium $\epsilon$, even if we tune the weight parameters $\lambda_g$ and $\lambda_d$ depending on $\epsilon$.

\subsection{The macroscopic auxiliary asymptotic-preserving neural networks (MA-APNNs)}
The above argument reveals that the vanilla PINNs are unsuitable for dealing with the multiscale feature.
In this subsection, we will devise a new asymptotic-preserving loss function to tackle this problem.
As a preliminary, we derive a macroscopic auxiliary equation from \eqref{eq:LRTEs-a}.
Then we demonstrate that the proposed loss function incorporating the macroscopic auxiliary equation is applicable to solve the LTREs for arbitrary $\epsilon$ and the solution converges to the diffusion limit as $\epsilon \rightarrow 0$.
\subsubsection{The macroscopic auxiliary equation}
From \eqref{eq:LRTEs-a}, we see that
\begin{equation}\label{iterate-sentence}
f = \left \langle f \right \rangle - \frac{\epsilon^2}{\sigma}\partial_t f - \frac{\epsilon}{\sigma}\Omega\cdot \nabla_r f - \frac{\epsilon^2}{\sigma}\alpha f + \frac{\epsilon^2}{\sigma}G.
\end{equation}
Replacing the $f$ in $\frac{\epsilon^2}{\sigma}\partial_t f$ and $\frac{\epsilon}{\sigma}\Omega\cdot \nabla_r f$ by the right-hand side of \eqref{iterate-sentence}, and arranging the result by the order of $\epsilon$, we get
%
\[
\begin{aligned}
& f = \langle f \rangle - \epsilon \frac{1}{\sigma}\Omega\cdot\nabla_r \left \langle f \right \rangle - \epsilon^2\left(\frac{1}{\sigma}\partial_t\left \langle f \right \rangle - \frac{1}{\sigma}\Omega\cdot\nabla_r(\frac{1}{\sigma}\Omega\cdot\nabla_r f)+\frac{1}{\sigma}\alpha f- \frac{1}{\sigma}G\right) \\
& +\epsilon^3\left(\frac{1}{\sigma}\partial_t(\frac{1}{\sigma}\Omega\cdot\nabla_rf) + \frac{1}{\sigma}\Omega\cdot\nabla_r \big(\frac{1}{\sigma}\partial_tf + \frac{1}{\sigma}\alpha f - \frac{1}{\sigma}G\big) \right)  \\
& +\epsilon^4\left(\frac{1}{\sigma}\partial_t \big( \frac{1}{\sigma}\partial_t f + \frac{1}{\sigma}\alpha f - \frac{1}{\sigma}G)\right).
\end{aligned}
\]
%
Let us pay attention to the second term with order $\epsilon^2$, namely, $\frac{\epsilon^2}{\sigma}\Omega\cdot\nabla_r(\frac{1}{\sigma}\Omega\cdot\nabla_r f)$.
We substitute the right side of \eqref{iterate-sentence} for the $f$ in this term, arrange the new equation by the order of $\epsilon$ and obtain
\begin{equation}\label{iterate-sentence2}
\begin{aligned}
f=&\left \langle f \right \rangle - \frac{\epsilon}{\sigma}\Omega\cdot\nabla_r \langle f\rangle - \frac{\epsilon^2}{\sigma} \left(\partial_t \langle f  \rangle - \Omega\cdot\nabla_r(\frac{1}{\sigma}\Omega\cdot\nabla_r \langle f \rangle) + \alpha f - G\right) \\
&+ \epsilon^3 \mathcal{A}(f,G) + \epsilon^4 \mathcal{B}(f,G),
\end{aligned}
\end{equation}
where
\[
\begin{aligned}
\mathcal{A}(f,G) &=\frac{1}{\sigma}\partial_t \big(\frac{1}{\sigma}\Omega\cdot\nabla_r f \big) + \frac{1}{\sigma} \Omega\cdot\nabla_r \left( \frac{1}{\sigma} \Big( \partial_tf + \alpha f - G -  \Omega\cdot\nabla_r(\frac{1}{\sigma}\Omega\cdot\nabla_r f)\Big) \right),
\end{aligned}
\]
\[
\begin{aligned}
\mathcal{B}(f,G)&=\frac{1}{\sigma^2} \partial_t \big(\partial_t f + \alpha f - G \big)-\frac{1}{\sigma}\Omega\cdot\nabla_r \left( \frac{1}{\sigma}\Omega\cdot\nabla_r \Big( \frac{1}{\sigma} \big( \partial_t f + \alpha f - G \big) \Big) \right).
\end{aligned}
\]
Integrating \eqref{iterate-sentence2} with respect to $\Omega$, in view of $\rho(t,r) = \langle f \rangle$ and $\langle \rho \rangle = \rho$, we arrive
\begin{equation}\label{eq:marco-auxiliary}
\partial_t\rho - \langle \Omega^2 \rangle \nabla_r\left(\frac{1}{\sigma}\nabla_r\rho\right) + \alpha\rho - G -\epsilon \langle \sigma \mathcal{A}(f,G) \rangle - \epsilon^2 \langle \sigma \mathcal{B}(f,G) \rangle =0,
\end{equation}
where $\langle \Omega^2  \rangle := \langle \mu^2 \rangle$ and we have used
\[
\begin{aligned}
& \langle \mu \rangle = \langle \xi \rangle = \langle \eta \rangle = 0 \quad (\text{i.e., } \langle \Omega \rangle = 0), \\
& \langle \mu^2 \rangle = \langle \xi^2 \rangle = \langle \eta^2 \rangle, \quad  \langle \mu \xi  \rangle = \langle \mu \eta  \rangle = \langle \xi \eta  \rangle =0.
\end{aligned}
\]
We call \eqref{eq:marco-auxiliary} the macroscopic auxiliary equation, which holds true if the solution of \eqref{eq:LRTEs-a} is smooth enough such that the high-order differentials in  \eqref{eq:marco-auxiliary} make sense.
When $\epsilon$ is tiny, we can regard \eqref{eq:marco-auxiliary} as the diffusion limit equation \eqref{eq:Diff-lim} with the perturbation $\epsilon \langle \sigma \mathcal{A}(f,G) \rangle + \epsilon^2 \langle \sigma \mathcal{B}(f,G) \rangle$.
When $\epsilon \rightarrow 0$, a formal calculation indicates that \eqref{eq:LRTEs-a} and \eqref{eq:marco-auxiliary} approaches the expected diffusion limit system:
\[
\sigma \left(\rho-f\right)=0, \quad \partial_t\rho - \left \langle \Omega^2 \right\rangle \nabla_r\left(\frac{1}{\sigma}\nabla_r\rho\right) + \alpha\rho -G =0.
\]
A rigorous validation of the convergence is out of the scope of this work. Our motivation is to apply the macroscopic auxiliary equation to design an asymptotic-preserving loss function.
\subsection{The exponentially weighted MA-APNNs loss function}
We use a deep neural network $f^{nn}_{\theta}(t,r,\Omega)$ to parameterize the radiation intensity $f(t,r,\Omega)$. In particular,  we choose the activation function $\sigma^o(\bm{x})=e^{-\bm{x}}$ at the output layer to guarantee the nonnegativity.
Now we design a new APNN loss $L_{\text{MA-APNNs}}^{\epsilon}$.
\begin{equation}\label{eq:MAAPNNs-loss}
L_{\text{MA-APNNs}}^{\epsilon} = L_{\text{MA-APNNs},g}^{\epsilon} + L_{\text{MA-APNNs},i}^{\epsilon} + L_{\text{MA-APNNs},b}^{\epsilon} + L_{\text{PINNs},c}^{\epsilon},
\end{equation}
where the residual of the mass conservation law, i.e., $L_{\text{PINNs},c}^{\epsilon}$, is the same to that of Section~\ref{eq:subsec-PINNs}.
The differences to the vanilla PINNs are shown as follows.

We add boundary and initial value constraints of the incident radiation $\left \langle f \right \rangle = \rho$ to $L_{\text{PINNs},b}^{\epsilon}$ and $L_{\text{PINNs},i}^{\epsilon}$, respectively.
\[
\begin{aligned}
& L_{\text{MA-APNNs},b}^{\epsilon}  = \lambda_{b} \left( \left\| B f_{\theta}^{nn} - f_b \right\|^2_{L^2(\tau \times \partial D \times S^2)} + \left\| B \left \langle f_{\theta}^{nn} \right \rangle -\left \langle f_b \right \rangle \right\|^2_{L^2(\tau \times \partial D)} \right), \\
& L_{\text{MA-APNNs},i}^{\epsilon}  = \lambda_{i} \left( \left\| f_{\theta}^{nn}(0) - f_0 \right\|^2_{L^2(D \times S^2)} + \left\| \left \langle f_{\theta}^{nn}(0) \right \rangle - \left \langle f_0 \right \rangle \right\|^2_{L^2(D)} \right).
\end{aligned}
\]

Most importantly, we replace the residual of the governing equation $L_{\text{PINNs},g}$ by
\[
\begin{aligned}
L_{\text{MA-APNNs},g}^{\epsilon} & = \Big \Vert \lambda_{\nu, \beta}^\frac{1}{2} \big( \epsilon^2\partial_t f_{\theta}^{nn} + \epsilon \Omega\cdot \nabla_r f_{\theta}^{nn}- \sigma
\left(\left \langle f_{\theta}^{nn} \right \rangle - f_{\theta}^{nn} \right) + \epsilon^2\alpha f_{\theta}^{nn} - \epsilon^2 G \big) \Big \Vert^2_{L^2(\tau \times D \times S^2)}  \\
& \quad + \Big \Vert (1-\lambda_{\nu,\beta})^\frac{1}{2} \Big( \partial_t \left \langle f_{\theta}^{nn} \right \rangle - \left\langle \Omega^2 \right \rangle \nabla_r \big( \frac{1}{\sigma}\nabla_r \left \langle f_{\theta}^{nn} \right \rangle \big) + \alpha \left \langle f_{\theta}^{nn} \right \rangle - G \\
& \qquad - \epsilon \left \langle \mathcal{A}(f_{\theta}^{nn},G) \right \rangle - \epsilon^2 \left \langle \mathcal{B}(f_{\theta}^{nn},G) \right \rangle \Big) \Big \Vert^2_{L^2(\tau \times D)},
\end{aligned}
\]
where $\lambda_{\nu,\beta}(r) := e^{-\nu(r)\beta_1} + \beta_2$ with $\nu(r) = \frac{\sigma(r)}{\epsilon^2}+\alpha(r)$. Here, $\beta=(\beta_1,\beta_2) > 0$ is the tunable parameter satisfying $1- \max_{r}\lambda_{\nu,\beta}(r) > 0$.

\begin{figure}[!htbp]
	\setlength{\abovecaptionskip}{0.cm}
	\setlength{\belowcaptionskip}{-0.cm}
	\centering
	{
		\begin{minipage}{4.0in}
			\centering
			\includegraphics[width=4.0in]{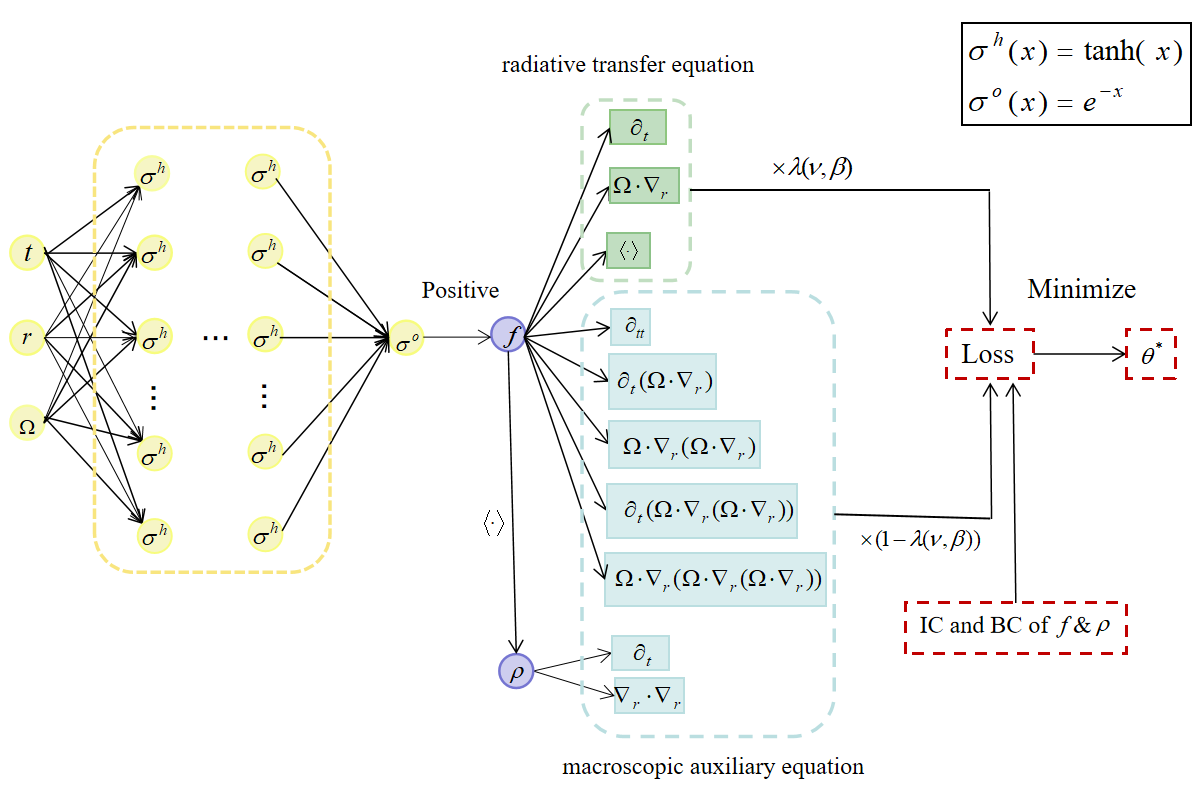}
		\end{minipage}
	}
	\caption{Schematic of MA-APNNs for solving the linear radiative transfer equation with initial and boundary data.}
	\label{figure3.1}
\end{figure}

It is easy to verify that the loss $L_{\text{MA-APNNs},g}^{\epsilon}$ asymptotically converges to the loss of the diffusion limit equations.
In fact, since $\sigma(r) > 0$, as $\epsilon \rightarrow 0$, we observe that $\lambda_{\nu, \beta}(r) \rightarrow \beta_2$ and
\[
\begin{aligned}
L_{\text{MA-APNNs},g}^\epsilon & \rightarrow \beta_2\left\| - \sigma
\left(\left \langle f_{\theta}^{nn} \right \rangle - f_{\theta}^{nn} \right) \right\|^2_{L^2(\tau \times D \times S^2)} \\
& + \left(1-\beta_2\right)\left\| \partial_t \left \langle f_{\theta}^{nn} \right \rangle - \left\langle \Omega^2 \right \rangle \nabla_r\left(\frac{1}{\sigma}\nabla_r \left \langle f_{\theta}^{nn} \right \rangle \right) + \alpha \left \langle f_{\theta}^{nn} \right \rangle - G \right\|^2_{L^2(\tau \times D)},
\end{aligned}
\]
which is the residual of the governing equation of the desired diffusion limit system \eqref{eq:Diff-lim-all}.
Since the smooth solution $f$ satisfies the macroscopic auxiliary equation \eqref{eq:marco-auxiliary}, the loss \eqref{eq:MAAPNNs-loss} has the consistency, that is, $L_{\text{MA-APNNs}}^{\epsilon}|_{f_{\theta}^{nn} =f} = 0$. Therefore, the loss function $L_{\text{MA-APNNs}}^\epsilon$ possesses the AP property and is applicable to the cases under the transport, medium, and diffusion regimes.
We put a schematic plot of our method in Figure~\ref{figure3.1}.

\subsection{The MA-APNN empirical loss function}
In practice, we apply the quadrature rules to compute the integrals of the loss function $L_{\text{MA-APNNs}}^{\epsilon}$.
\subsubsection{Angular integration approximation}
We utilize the Gauss-Legendre quadrature rule to approximate the angular integration $\langle f \rangle$:
\[
\left \langle f \right \rangle=\frac{1}{4\pi}\int_{S^2} f(t,r,\Omega) \text{d}\Omega \approx \frac{1}{4\pi} \sum_{m=1}^{N_s} f(\Omega_m) \omega_m.
\]
The Gauss-Legendre quadrature points $\{\Omega_m\}_{m=1}^{N_s} \subset S^2$ will be used as the training points as well, and $\{\omega_m\}_{m=1}^{N_s}$ are the associated weights.
\subsubsection{Empirical loss functions, sample points, and labeled data}
Here we sample the training points $\wp_{f}=\left\{\bm{x}^f_j=(t^f_j, r^f_j, \Omega^f_j), j=1,2,...N_f\right\}$ as a low-discrepancy Sobol sequence in the computational domain, and choose the Quasi-Monte Carlo method for numerical integration \cite{sobol1967distribution}.
The empirical loss function of the governing equation is stated as follows
\[
\begin{aligned}
& L_{\text{MA-APNNs},g}^{\epsilon,nn} =
\frac{1}{N_{\text{int}}}\sum_{j=1}^{N_{\text{int}}} \Big| \lambda_{\nu,\beta}(r_j^{\text{int}}) \Big (\epsilon^2\partial_t f_{\theta}^{nn}(t_j^{\text{int}}, r_j^{\text{int}}, \Omega_j^{\text{int}}) + \epsilon \Omega\cdot \nabla_r f_{\theta}^{nn}(t_j^{\text{int}}, r_j^{\text{int}}, \Omega_j^{\text{int}}) \\
&-\sigma \left( \langle f_{\theta}^{nn}\rangle(t_j^{\text{int}}, r_j^{\text{int}}) - f_{\theta}^{nn}(t_j^{\text{int}}, r_j^{\text{int}}, \Omega_j^{\text{int}}) \right) + \epsilon^2\alpha f_{\theta}^{nn}(t_j^{\text{int}}, r_j^{\text{int}}, \Omega_j^{\text{int}}) - \epsilon^2 G \Big ) \Big|^2  \\
&+ \frac{1}{N_{\text{int}}}\sum_{j=1}^{N_{\text{int}}} \Big| \left(1-\lambda_{\nu,\beta}(r_j^{\text{int}}) \right) \Big( \partial_t \langle f_{\theta}^{nn}\rangle(t_j^{\text{int}}, r_j^{\text{int}})  - \langle \Omega^2 \rangle \nabla_r\left(\frac{1}{\sigma}\nabla_r \langle f_{\theta}^{nn}\rangle(t_j^{\text{int}}, r_j^{\text{int}}) \right) \\
&+ \alpha \langle f_{\theta}^{nn}\rangle(t_j^{\text{int}}, r_j^{\text{int}})  -G -\epsilon \langle \sigma \mathcal{A}(f_{\theta}^{nn},G) \rangle (t_j^{\text{int}}, r_j^{\text{int}}) - \epsilon^2 \langle \sigma \mathcal{B}(f_{\theta}^{nn},G) \rangle(t_j^{\text{int}}, r_j^{\text{int}}) \Big ) \Big|^2
\end{aligned}
\]
where $\{(t_j^{\text{int}},r_j^{\text{int}},\Omega_j^{\text{int}})\}_{j=1}^{N_{\text{int}}}$ are interior Sobol sequence points.
Similarly, the empirical loss functions for the boundary condition, the initial value, and conservation law are given by:
\[
\begin{aligned}
L_{\text{MA-APNNs},b}^{\epsilon,nn} = \lambda_{b} \frac{1}{N_{\text{sb}}}\sum_{j=1}^{N_{\text{sb}}} &\Big( \left| B f_{\theta}^{nn}(t_j^{\text{sb}}, r_j^{\text{sb}}, \Omega_j^{\text{sb}})-f_b(t_j^{\text{sb}}, r_j^{\text{sb}}, \Omega_j^{\text{sb}}) \right|^2 \\
&+ \left| B \langle f_{\theta}^{nn} \rangle (t_j^{\text{sb}}, r_j^{\text{sb}}) - \langle f_b \rangle (t_j^{\text{sb}}, r_j^{\text{sb}}) \right|^2 \Big),
\end{aligned}
\]
\[
L_{\text{MA-APNNs},i}^{\epsilon,nn} = \lambda_i \frac{1}{N_{\text{tb}}}\sum_{j=1}^{N_{\text{tb}}} \left( \left| f_{\theta}^{nn}(0, r_j^{\text{tb}}, \Omega_j^{\text{tb}})-f_0(r_j^{\text{tb}}, \Omega_j^{\text{tb}}) \right|^2 + \left| \langle f_{\theta}^{nn} \rangle (0, r_j^{\text{tb}})- \langle f_0 \rangle (r_j^{\text{tb}})\right|^2 \right),
\]
\[
L_{\text{PINNs},c}^{\epsilon,nn} = \lambda_{c} \frac{1}{N_{\text{c}}}\sum_{j=1}^{N_{\text{c}}} \left| \epsilon \partial_t \langle f_{\theta}^{nn} \rangle(t_j^{\text{c}},r_j^{\text{c}})  + \langle \Omega \cdot n_r f_{\theta}^{nn} \rangle (t_j^{\text{c}},r_j^{\text{c}})  + \epsilon \alpha \langle f_{\theta}^{nn}\rangle (t_j^{\text{c}},r_j^{\text{c}}) - \epsilon G \right|^2,
\]
where $\{ (t_j^{\text{sb}},r_j^{\text{sb}},\Omega_j^{\text{sb}}) \}_{j=1}^{N_{\text{sb}}}$ are Sobol sequence points on boundary, $\{(0,r_j^{\text{tb}},\Omega_j^{\text{tb}})\}_{j=1}^{N_{\text{tb}}}$ the Sobol sequence points at $t=0$, and $\{ (t_j^{\text{c}},r_j^{\text{c}}) \}_{j=1}^{N_{\text{c}}}$ the Sobol sequence points on $\tau \times D$.

We define the total empirical loss function for MA-APNNs:
\begin{equation}
L_{\text{MA-APNNs}}^{\epsilon,nn} = L_{\text{MA-APNNs},g}^{\epsilon,nn} + L_{\text{MA-APNNs},b}^{\epsilon,nn} + L_{\text{MA-APNNs},i}^{\epsilon,nn}  + L_{\text{PINNs},c}^{\epsilon,nn}.
\end{equation}
Now, we are in the position to find the solution of the minimization problem
\begin{equation}
\theta^{*}= \mathop{\arg\min}\limits_{\theta}\left( L_{\text{MA-APNNs}}^{\epsilon,\text{nn}}\right).
\end{equation}
In the simulation, we adopt appropriate optimization algorithms (e.g. Adam, LBFGS) to minimize the non-convex loss function $L_{\text{MA-APNNs}}^{\epsilon,\text{nn}}$.
\section{Numerical Experiments}\label{sec4}
We present several examples to compare the performance of PINNs, APNNs based on micro-macro decomposition\cite{jin2023asymptotic1} and MA-APNNs for solving multiscale LRTEs.
The result confirms that the MA-APNNs have an advantage in solving diffusive scaling problems.
In the numerical simulations, we adopt the Adam algorithm with an initial learning rate $0.001$ in the optimization process of the loss function and use $\tanh x$ as the activation function of all hidden layers. Unless stated otherwise, the reference solution $f^{\text{Ref}}$ is obtained by the UGKS\cite{mieussens2013asymptotic}.
We calculate $L^2$-norm relative errors of $\rho= \langle f \rangle$:
\[
L^2_{\text{error}}(\rho) = \frac{\Vert \langle f^{nn} \rangle-\langle f^{\text{Ref}}\rangle \Vert_{L^{2}(\tau \times D )}}{\Vert \langle f^{\text{Ref}} \rangle \Vert_{L^{2}(\tau \times D)}}.
\]
\subsection{\textbf{One-dimensional problems}}
\subsubsection{\textbf{Kinetic regime with isotropic boundary condition}}\label{ex:1}
\begin{figure}[!htbp]
	\setlength{\abovecaptionskip}{0.cm}
	\setlength{\belowcaptionskip}{-0.cm}
	\centering
	{
		\begin{minipage}{1.5in}
			\centering
			\includegraphics[width=1.5in]{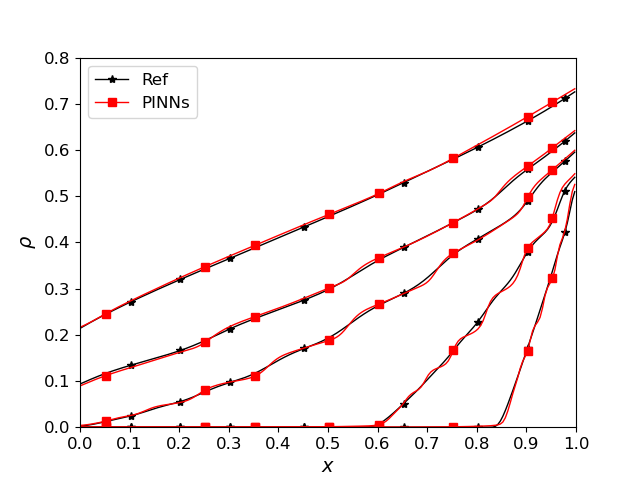}
		\end{minipage}
	}
	{
	\begin{minipage}{1.5in}
		\centering
		\includegraphics[width=1.5in]{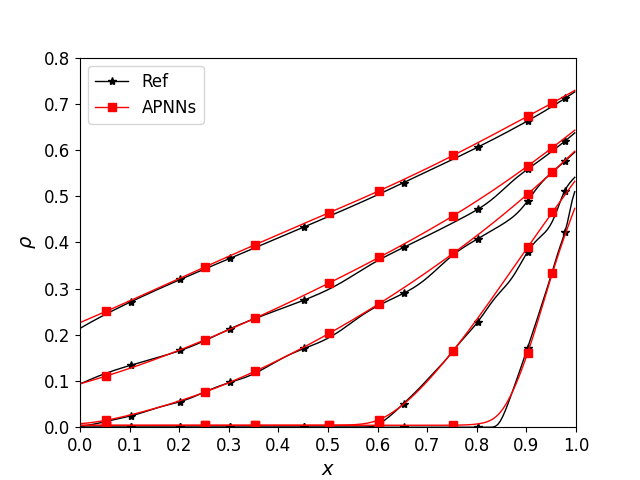}
	\end{minipage}
	}
	{
		\begin{minipage}{1.5in}
			\centering
			\includegraphics[width=1.5in]{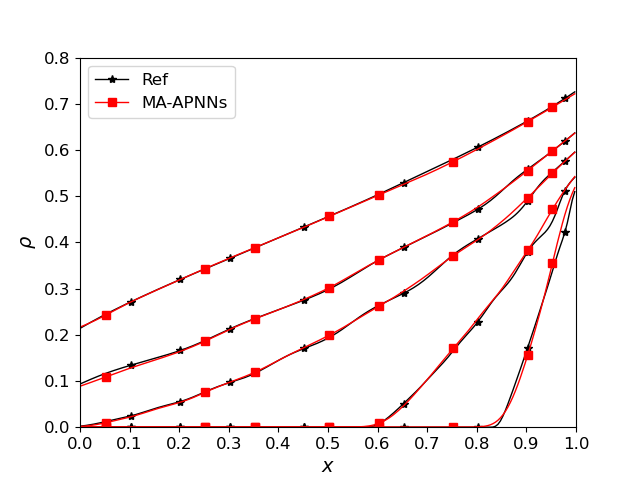}
		\end{minipage}
	}
	\caption{Kinetic regime with $\epsilon=1$. The density $\rho$ at times $t=0.15, 0.4, 1.0, 1.6, 4.0$. (Left) Ref v.s. PINNs. (Middle) Ref v.s. APNNs. (Right) Ref v.s. MA-APNNs.}
	\label{figure4.1.1}
\end{figure}

Consider the 1D LRTEs in the kinetic regime with $\epsilon = 1$.
We set $D= (0,1)$, $\tau = (0,4)$, $\sigma=1$, $\alpha = G=0$, and enforce the isotropic inflow boundary condition \eqref{eq:LRTEs-1d-inflow-bc} and the initial condition \eqref{eq:LRTEs-c} with
\[
f_L(t,\mu>0)=0, \quad  f_R(t,\mu<0)=1, \quad f_0(x,\mu)=0.
\]
We take the network structure $f_{\theta}^{nn}=[3,40,40,40,40,1]$ with the hyperparameters $\beta_1=10^{-3}$, $\beta_2=10^{-4}$, $\lambda_{b}=\lambda_i =1$, $\lambda_{c}=0$. The number of the batch training samples is $(N_{\text{int}},N_{\text{sb}},N_{\text{tb}})=(2000,400 \times 2,400)$. According to the content introduced in  \cite{jin2023asymptotic1}, we set APNNs here as $g_{\theta_1}^{nn}=[3,40,40,40,40,1]$, $\rho_{\theta_2}^{nn}=[2,40,40,40,40,1]$. The reference solution and the prediction obtained by PINNs, APNNs and MA-APNNs are plotted in Figure~\ref{figure4.1.1}. The $L^2$ errors and training time of PINNs, APNNs and MA-APNNs are displayed in  Table~\ref{Table4.1.1}. All three methods are well applicable to the LRTEs in the kinetic regime.
\begin{table}[!htbp]
	\centering
	\caption{Kinetic regime with $\epsilon=1$. The errors and training time of PINNs, APNNs and MA-APNNs.}
	\label{Table4.1.1}
	\scalebox{0.9}{
	\begin{tabular}{lllllll}
		\hline\noalign{\smallskip}
		$L^2_{\text{error}}(\rho)$ &$t=0.15$ & $t=0.4$ & $t=1.0$ & $t=1.6$ & $t=4.0$ & Training time \\
		\noalign{\smallskip}\hline\noalign{\smallskip}
		PINNs &3.48e-02   &2.86e-02   &1.41e-02  &1.20e-02 & 1.01e-02 &31min \\
		\noalign{\smallskip}\hline\noalign{\smallskip}
		APNNs &6.19e-02   &3.89e-02   &2.66e-02  &2.49e-02 & 1.47e-02 &1h 18min \\
		\noalign{\smallskip}\hline\noalign{\smallskip}
		MA-APNNs &6.51e-02   &2.86e-02   &1.22e-02  &8.60e-03 & 5.06e-03 &2h 15min \\
		\noalign{\smallskip}\hline
	\end{tabular}}
\end{table}
\subsubsection{\textbf{Initial layer problem with periodic boundary condition}}\label{ex:2}
We adopt the same settings of the previous example, except that we set the temporal domain $\tau=(0,1)$, replace the inflow boundary condition by the periodic boundary condition \eqref{eq:LRTEs-1d-periodic-bc}, and set the initial value
\[
f_0(x,\mu)=\frac{1+\cos(4\pi x)}{\sqrt{2\pi}}e^{-\frac{\mu^2}{2}}.
\]
The solution of this 1D LRTEs is known to have an initial layer neat $t=0$. We utilize MA-APNNs with soft and hard constraints of the periodic boundary condition and APNNs with boundary hard constraints to solve this problem.

\begin{figure}[!htbp]
	\setlength{\abovecaptionskip}{0.cm}
	\setlength{\belowcaptionskip}{-0.cm}
	\centering
	{
	\begin{minipage}{1.5in}
		\centering
		\includegraphics[width=1.5in]{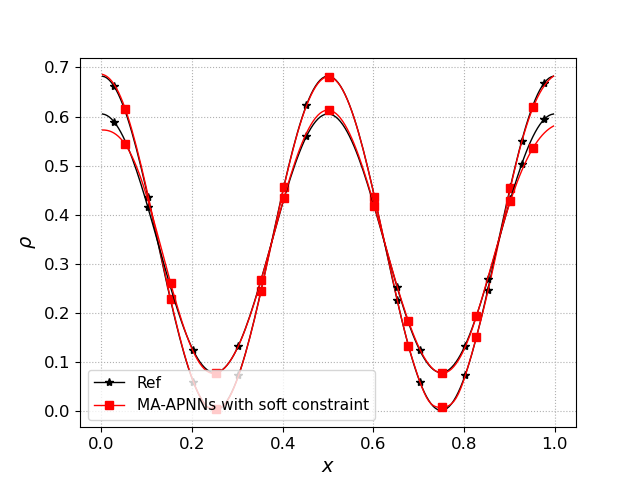} 
	\end{minipage}
	}
	{
	\begin{minipage}{1.5in}
		\centering
		\includegraphics[width=1.5in]{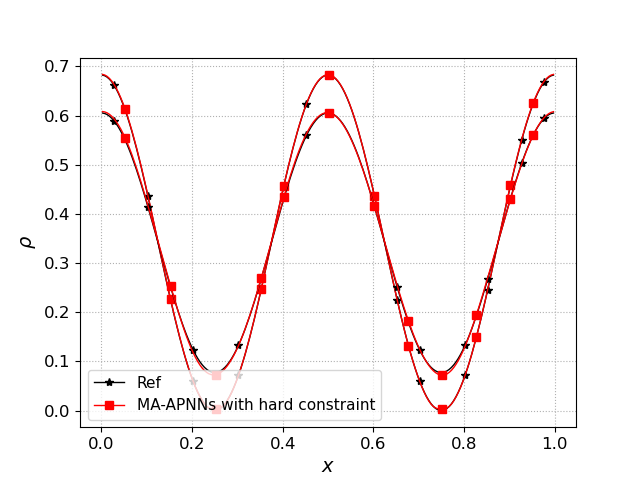} 
	\end{minipage}
	}
	{
		\begin{minipage}{1.5in}
			\centering
			\includegraphics[width=1.5in]{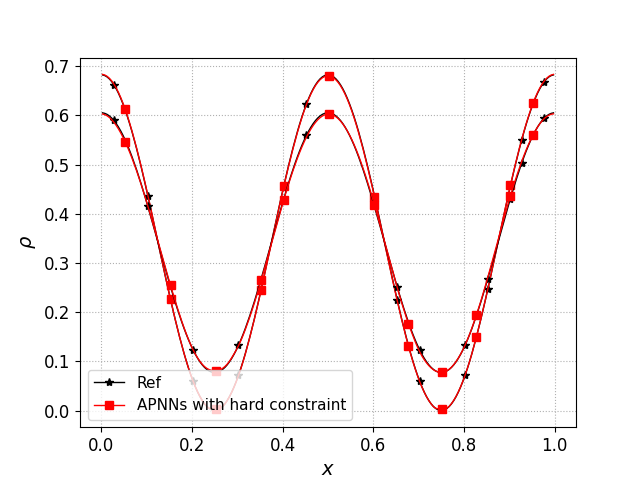}
		\end{minipage}
	}
	\caption{Initial layer with $\epsilon=1$. The density $\rho$ at times $t=0.0, 0.1$. (Left) Ref v.s. MA-APNNs with boundary soft constraint. (Middle) Ref v.s. MA-APNNs with boundary hard constraint. (Right) Ref v.s. APNNs with boundary hard constraint.}
	\label{figure4.1.2}
\end{figure}

We set $(N_{\text{int}},N_{\text{sb}},N_{\text{tb}})=(2000, 500\times 2, 1000)$, and $(\beta_1, \beta_2)=(10^{-3}, 10^{-5})$. Under the soft constraints of the boundary condition, we choose $f_{\theta}^{nn}=[3,40,40,40,40,1]$, $(\lambda_{b},\lambda_c,\lambda_i)=(1,1,1000)$. Under the hard constraints, we add a custom layer between the input and the first hidden layer in MA-APNNs and APNNs to transform $x$ to $(\sin 2\pi x, \cos 2\pi x)$\cite{mojgani2023kolmogorov,dong2021method} such that the periodic boundary condition is satisfied for $f_{\theta}^{nn}$ in MA-APNNs and $g_{\theta_1}^{nn}$, $\rho_{\theta_2}^{nn}$ in APNNs\cite{jin2023asymptotic1}. We take $f_{\theta}^{nn}/g_{\theta_1}^{nn}=[4,40,40,40,40,1]$, $\rho_{\theta_2}^{nn}=[3,40,40,40,40,1]$ and $(\lambda_{b},\lambda_c,\lambda_i)=(0,1,1000)$. The numerical solutions of APNNs and MA-APNNs are plotted in Figure~\ref{figure4.1.2}, together with the reference solution. The $L^2$ errors and training time are shown in Table~\ref{Table4.1.2}. We observe that the MA-APNNs with hard constraints of the boundary condition perform slightly better than the soft constraints. APNNs and MA-APNNs with hard constraints both perform well under this case.
\begin{table}[!htbp]
	\centering
	\caption{Initial layer with $\epsilon=1$. The errors and training time of MA-APNNs with boundary soft and hard constraint and APNNs with boundary hard constraint.}
	\label{Table4.1.2}
	\scalebox{0.9}{
		\begin{tabular}{llll}
			\hline\noalign{\smallskip}
			$L^2_{\text{error}}(\rho)$ &$t=0.0$ & $t=0.1$ &Training time \\
			\noalign{\smallskip}\hline\noalign{\smallskip}
			MA-APNNs with soft constraint &6.38e-03    &2.55e-02 &3h 10min\\
			\noalign{\smallskip}\hline\noalign{\smallskip}
			MA-APNNs with hard constraint &2.31e-03    & 7.31e-03 &3h 13min \\
			\noalign{\smallskip}\hline\noalign{\smallskip}
			APNNs with hard constraint &2.89e-03    &6.28e-03 &1h 29min\\
			\noalign{\smallskip}\hline
	\end{tabular}}
\end{table}

\subsubsection{\textbf{Diffusion regime with a constant scattering frequency}}\label{ex:3}
\begin{figure}[!htbp]
	\setlength{\abovecaptionskip}{0.cm}
	\setlength{\belowcaptionskip}{-0.cm}
	\centering
	{
	\begin{minipage}{1.5in}
		\centering
		\includegraphics[width=1.5in]{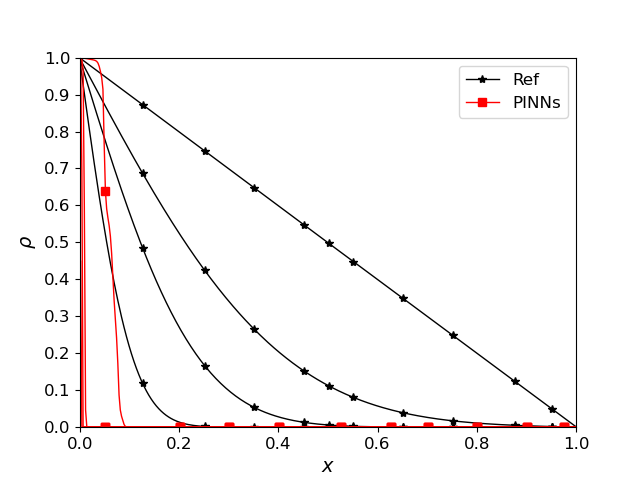}
	\end{minipage}
	}
	{
		\begin{minipage}{1.5in}
			\centering
			\includegraphics[width=1.5in]{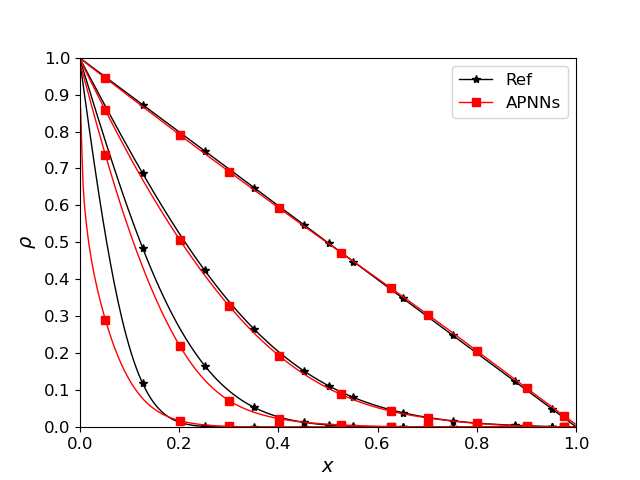}
		\end{minipage}
	}
	{
		\begin{minipage}{1.5in}
			\centering
			\includegraphics[width=1.5in]{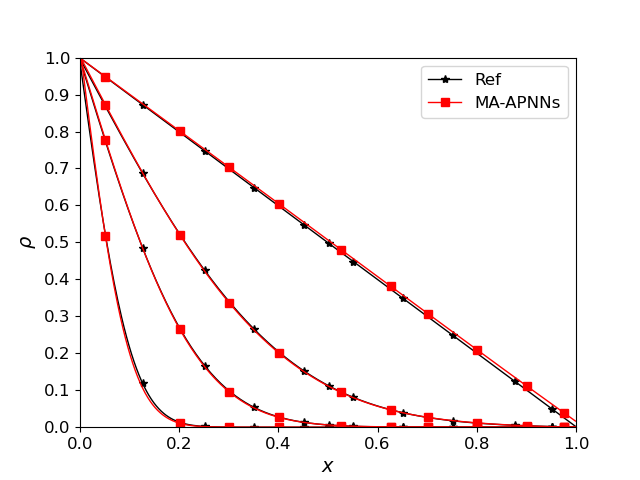}
		\end{minipage}
	}
	\caption{Diffusion regime with $\epsilon=10^{-8}$. The density $\rho$ at times $t=0.01, 0.05, 0.15, 2.00$. (Left) Ref v.s. PINNs. (Middle) Ref v.s. APNNs. (Right) Ref v.s. MA-APNNs.}
	\label{figure4.1.3}
\end{figure}

Consider the 1D LRTEs under the following settings:
\[
\begin{aligned}
& D = (0,1), \quad  \tau=(0,2), \quad  f_L(t,\mu>0)=1, \quad  f_R(t,\mu<0)=0,  \\
& f_0(x,\mu)=0, \quad \sigma=1, \quad \alpha=0, \quad G=0, \quad \epsilon=10^{-8}.
\end{aligned}
\]
Since the Knudsen number $\epsilon$ is extremely small (in diffusion scale), the solution is close to the diffusion limit of \eqref{eq:Diff-lim-all}.
We choose $f_{\theta}^{nn}=[3,40,40,40,40,1]$ in PINNs and MA-APNNs, $ g_{\theta_1}^{nn}=[3,40,40,40,40,1],\rho_{\theta_2}^{nn}=[2,40,40,40,40,1]$ in APNNs\cite{jin2023asymptotic1} and set $(N_{\text{int}},N_{\text{sb}},N_{\text{tb}})=(1000, 200\times 2, 200)$ and $(\beta_1,\beta_2, \lambda_{b}, \lambda_i, \lambda_c)=(10^{-5}, 10^{-16}, 10, 1, 0)$.
The prediction of PINNs, APNNs and MA-APNNs are presented in Figure~\ref{figure4.1.3}. The $L^2$ errors and training time are listed in Table~\ref{Table4.1.3}.
The experimental result indicates that the vanilla PINNs fail to approximate the diffusion scale problem, and APNNs do badly near initial times while the proposed MA-APNN method is well applicable to the case with very tiny $\epsilon$.
\begin{table}[!htbp]
	\centering
	\caption{Diffusion regime with $\epsilon=10^{-8}$. The errors and training time of PINNs, APNNs and MA-APNNs.}
	\label{Table4.1.3}
	\scalebox{0.9}{
		\begin{tabular}{llllll}
			\hline\noalign{\smallskip}
			$L^2_{\text{error}}(\rho)$ &$t=0.01$ & $t=0.05$ & $t=0.15$ & $t=2.00$ & Training time \\
			\noalign{\smallskip}\hline\noalign{\smallskip}
			PINNs &9.75e-01  &9.71e-01  &9.65e-01  &8.96e-01 &17min 34s\\
			\noalign{\smallskip}\hline\noalign{\smallskip}
			APNNs &3.81e-01  &8.18e-02  &2.37e-02  &8.95e-03 &2h 45min \\
			\noalign{\smallskip}\hline\noalign{\smallskip}
			MA-APNNs &3.16e-02  &5.79e-03  &5.30e-03  &1.35e-02 &1h 24min \\
			\noalign{\smallskip}\hline
	\end{tabular}}
\end{table}
\subsubsection{\textbf{Diffusion regime with a variable scattering frequency}}\label{ex:4}
\begin{figure}[!htbp]
	\centering
	{
	\begin{minipage}{1.5in}
		\centering
		\includegraphics[width=1.5in]{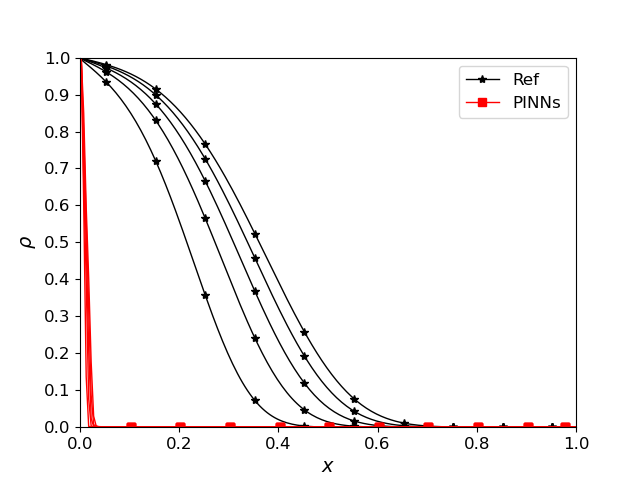}
	\end{minipage}
	}
	{
	\begin{minipage}{1.5in}
		\centering
		\includegraphics[width=1.5in]{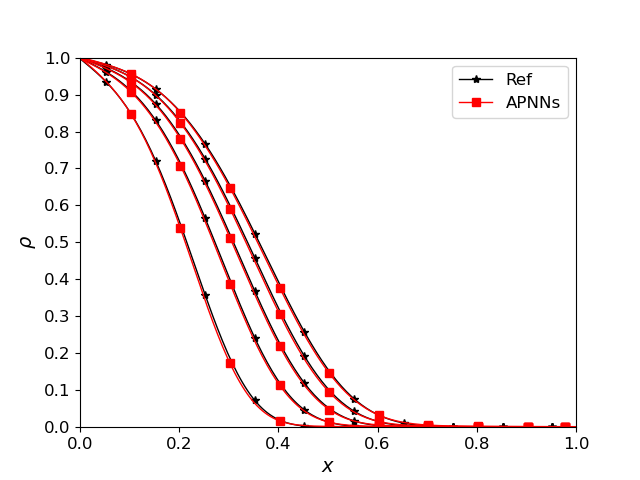}
	\end{minipage}
}
	{
		\begin{minipage}{1.5in}
			\centering
			\includegraphics[width=1.5in]{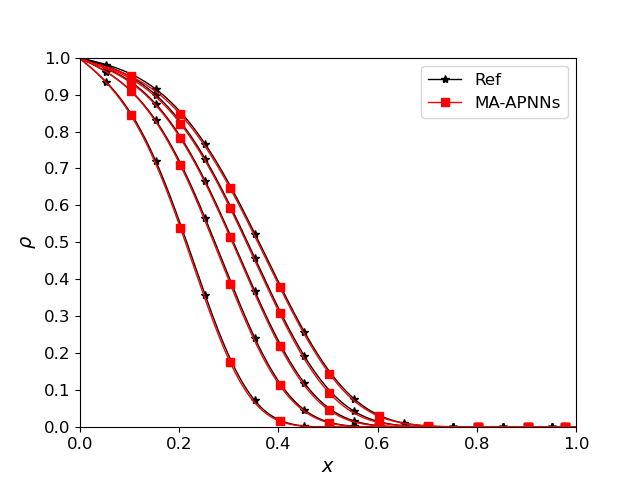}
		\end{minipage}
	}
	\caption{Diffusion regime with $\epsilon=10^{-4}$. The density $\rho$ at times $t=0.2, 0.4, 0.6, 0.8, 1.0$. (Left) Ref v.s. PINNs. (Middle) Ref v.s. APNNs. (Right) Ref v.s. MA-APNNs.}
	\label{figure4.1.4}
\end{figure}

Use the same settings as the previous example except that we set the scattering frequency $\sigma=1+(10x)^2$ and the Knudsen number $\epsilon=10^{-4}$.
In numerical simulation, the network structures are the same as in \ref{ex:3} and we take $(N_{\text{int}},N_{\text{sb}},N_{\text{tb}})=(1000, 200\times 2, 400)$, $(\beta_1,\beta_2, \lambda_{b}, \lambda_i, \lambda_c)=(10^{-5}, 10^{-16}, 1, 1, 0)$.
We see that the prediction of APNNs and MA-APNNs coincide with the reference solutions (see Figure~\ref{figure4.1.4} (middle and right)) while PINNs fail the test (see Figure~\ref{figure4.1.4} (left)).
The $L^2$ errors and training time are shown in Table~\ref{Table4.1.4}. MA-APNNs have less training time with comparable prediction accuracy to APNNs\cite{jin2023asymptotic1}.
\begin{table}[!htbp]
	\centering
	\caption{Diffusion regime with $\epsilon=10^{-4}$. The errors and training time of PINNs, APNNs and MA-APNNs.}
	\label{Table4.1.4}
	\scalebox{0.93}{
	\begin{tabular}{lllllll}
		\hline\noalign{\smallskip}
		$L_{\text{error}}^2(\rho)$ &$t=0.2$ & $t=0.4$ & $t=0.6$ & $t=0.8$ & $t=1.0$ &Training time\\
		\noalign{\smallskip}\hline\noalign{\smallskip}
		PINNs &9.65e-01    &9.66e-01     &9.66e-01   &9.65e-01 &9.64e-01  &28min 36s\\
		\noalign{\smallskip}\hline\noalign{\smallskip}
		APNNs &1.46e-02    &1.01e-02     &9.18e-03   &9.40e-03 &8.06e-03  &4h 26min \\
		\noalign{\smallskip}\hline\noalign{\smallskip}
		MA-APNNs &1.28e-02    &9.42e-03     &8.25e-03   &9.49e-03 &9.76e-03  &2h 25min\\
		\noalign{\smallskip}\hline
	\end{tabular}}
\end{table}
\subsubsection{\textbf{Intermediate regime with a variable scattering frequency and source term}}\label{ex:5}
\begin{figure}[!htbp]
	\centering
	{
	\begin{minipage}{1.5in}
		\centering
		\includegraphics[width=1.5in]{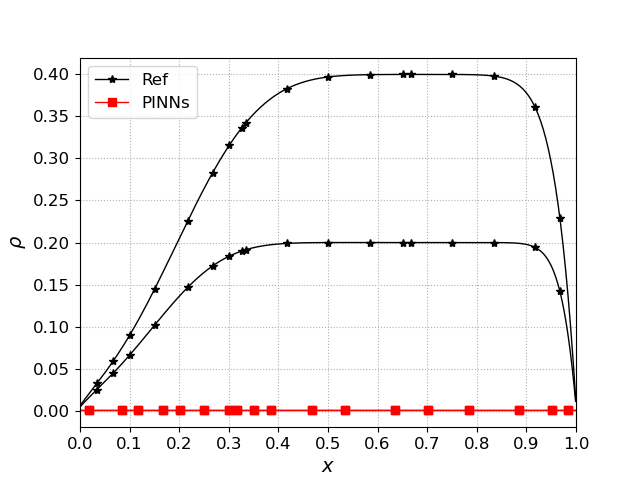}
	\end{minipage}
	}
	{
	\begin{minipage}{1.5in}
		\centering
		\includegraphics[width=1.5in]{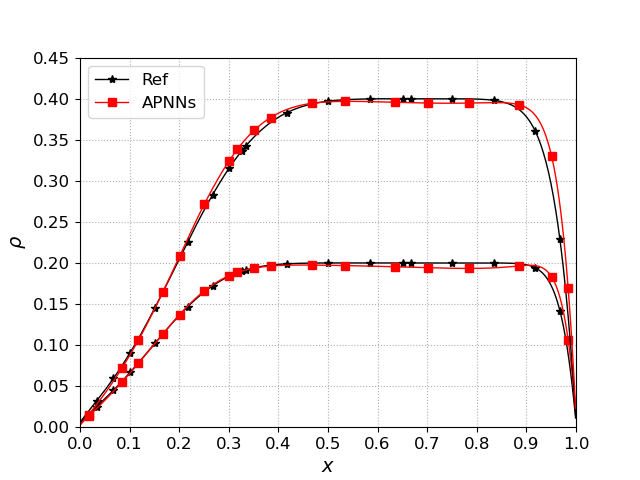}
	\end{minipage}
	}
	{
		\begin{minipage}{1.5in}
			\centering
			\includegraphics[width=1.5in]{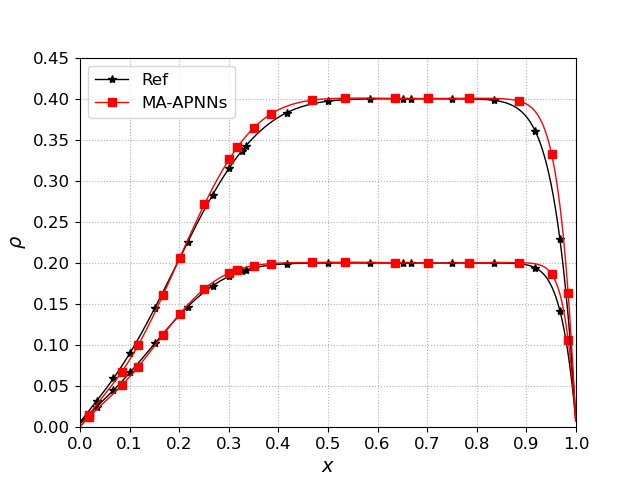}
		\end{minipage}
	}
	\caption{Intermediate regime with $\epsilon=10^{-2}$. The density $\rho$ at times $t=0.2, 0.4$. (Left) Ref v.s. PINNs. (Middle) Ref v.s. APNNs. (Right) Ref v.s. MA-APNNs.}
	\label{figure4.1.5}
\end{figure}

We replace the values of the source term and Knudsen number of the previous example by $G=1$ and $\epsilon=10^{-2}$, respectively.
We adopt the same settings of the network structure and hyperparameters except for $\beta_2=10^{-12}$.
We apply PINNs, APNNs and MA-APNNs to solve this problem. Specially, we set the activation functions of these networks' output layers as unit functions. The numerical results are shown in Figure~\ref{figure4.1.5} and Table~\ref{Table4.1.5}.
Together with the above examples, we are convinced that APNNs and MA-APNNs are able to resolve the multiscale characteristics of the 1D LRTEs, while the vanilla PINNs can only treat the case $\epsilon = O(1)$.
\begin{table}[!htbp]
	\centering
	\caption{Intermediate regime with $\epsilon=10^{-2}$. The errors and training time of PINNs, APNNs and MA-APNNs.}
	\label{Table4.1.5}
	\begin{tabular}{llll}
		\hline\noalign{\smallskip}
		$L_{\text{error}}^2(\rho)$ &$t=0.2$ & $t=0.4$ &Training time\\
		\noalign{\smallskip}\hline\noalign{\smallskip}
		PINNs &9.98e-01    &9.99e-01  &30min 9s   \\
		\noalign{\smallskip}\hline\noalign{\smallskip}
		APNNs &2.87e-02    &3.24e-02  &4h 37min  \\
		\noalign{\smallskip}\hline\noalign{\smallskip}
		MA-APNNs &2.70e-02    &3.44e-02 &2h 25min    \\
		\noalign{\smallskip}\hline
	\end{tabular}
\end{table}
\subsection{\textbf{Two-dimensional problems}}
We apply the MA-APNNs to solve the 2D LRTEs under the following settings
\[
D=(0,1)^2, \quad \tau = (0,1), \quad \sigma=1, \quad \alpha=0, \quad G=1.
\]
The homogeneous initial and inflow boundary conditions are imposed, i.e., $f_0 = 0$ and $f_B|_{\Gamma_-} = 0$.
\subsubsection{\textbf{Kinetic regime ($\epsilon=1$)}}\label{ex:7}
\begin{figure}[!htbp]
\begin{tabular}{ccc}
\includegraphics[width=0.32\textwidth]{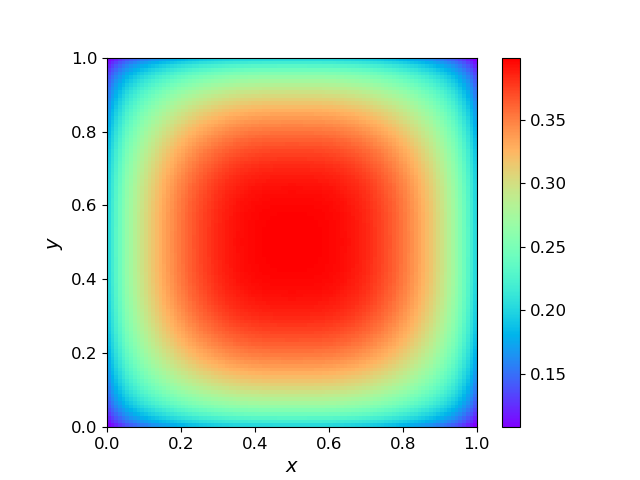} &
\includegraphics[width=0.32\textwidth]{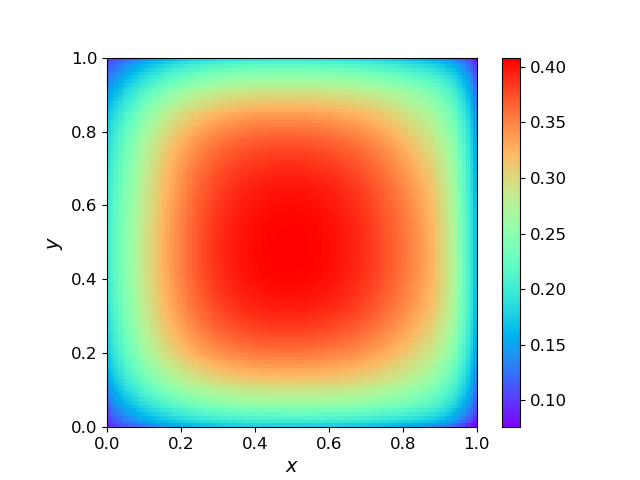} &
\includegraphics[width=0.32\textwidth]{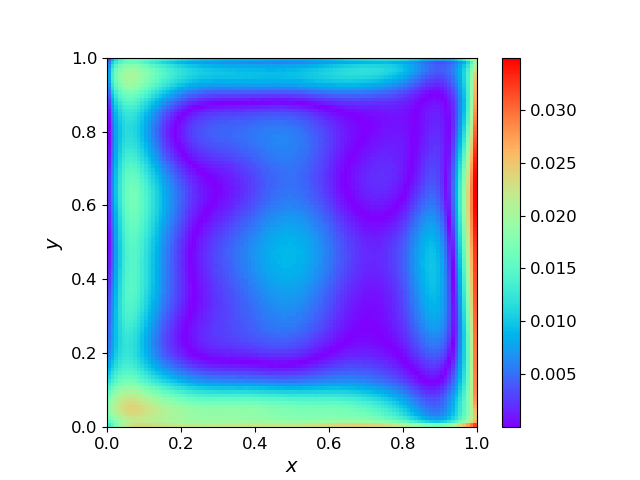} \\
(a) $\rho^{\text{Ref}}$, $t=0.4$ & (b) $\rho^{nn}_\theta$, $t=0.4$ & (c) $|\rho^{\text{Ref}}-\rho^{nn}_\theta|$, $t=0.4$ \\
\includegraphics[width=0.32\textwidth]{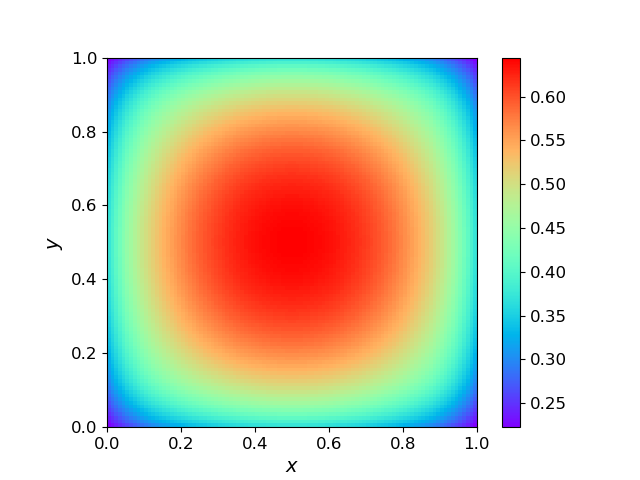} &
\includegraphics[width=0.32\textwidth]{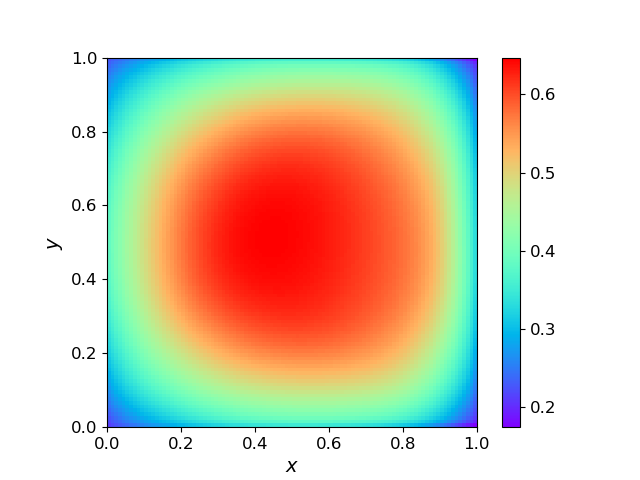} &
\includegraphics[width=0.32\textwidth]{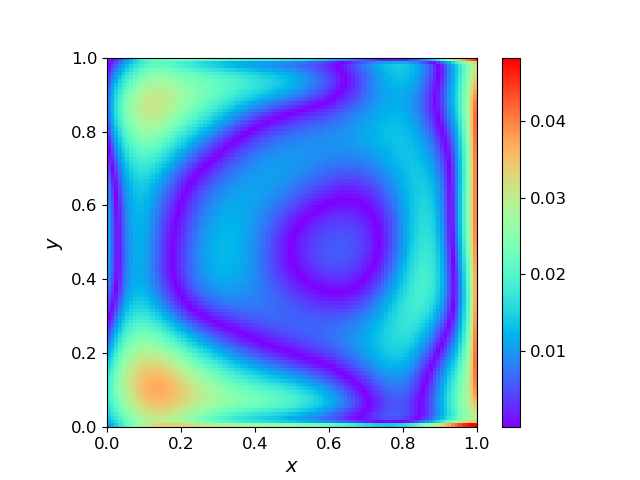} \\
(d) $\rho^{\text{Ref}}$, $t=1.0$ & (e) $\rho^{nn}_\theta$, $t=1.0$ & (f) $|\rho^{\text{Ref}}-\rho^{nn}_\theta|$, $t=1.0$ \\
\end{tabular}
\caption{The 2D LRTEs in the kinetic regime ($\epsilon=1$).}
\label{figure4.2.1}
\end{figure}

We set $\tilde{f}_{\theta}^{nn}=[5,40,40,40,40,1]$ and the boundary and initial conditions are enforced as hard constraints, that is,
\[
f_\theta^{nn}= t (x+ \text{Relu}(-\xi)^2)(1-x+\text{Relu}(\xi)^2) (y+\text{Relu}(-\eta)^2)(1-y+\text{Relu}(\eta)^2) \tilde{f}_{\theta}^{nn}.
\]
And we take $N_{\text{int}}=2000$ and $(\beta_1,\beta_2,\lambda_{b},\lambda_i,\lambda_{c})=(10^{-6}, 10^{-7},0,0,0)$.
The reference solution and the prediction obtained by MA-APNNs are plotted in Figure~\ref{figure4.2.1}.
The $L^2$ relative errors at times $t=0.4, 1.0$ are $2.91e-02$ and $2.84e-02$, respectively.
%

\subsubsection{\textbf{Diffusion regime} $(\epsilon=10^{-8})$}\label{ex:8}
\begin{figure}[!htbp]
\begin{tabular}{ccc}
\includegraphics[width=0.32\textwidth]{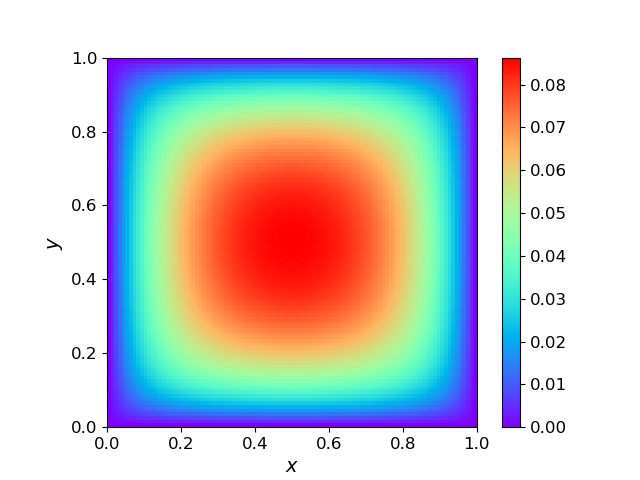} &
\includegraphics[width=0.32\textwidth]{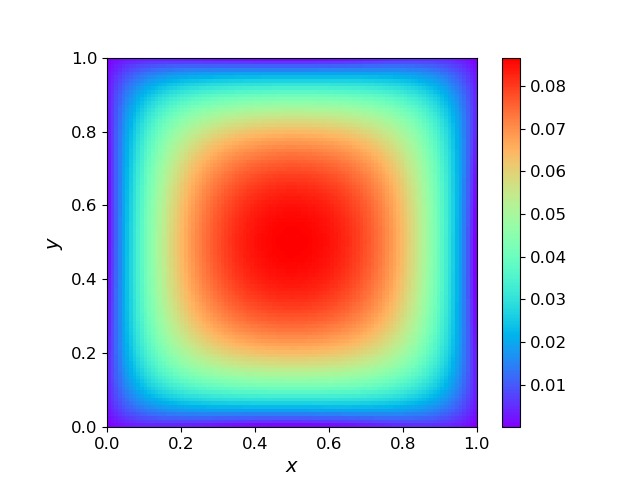} &
\includegraphics[width=0.32\textwidth]{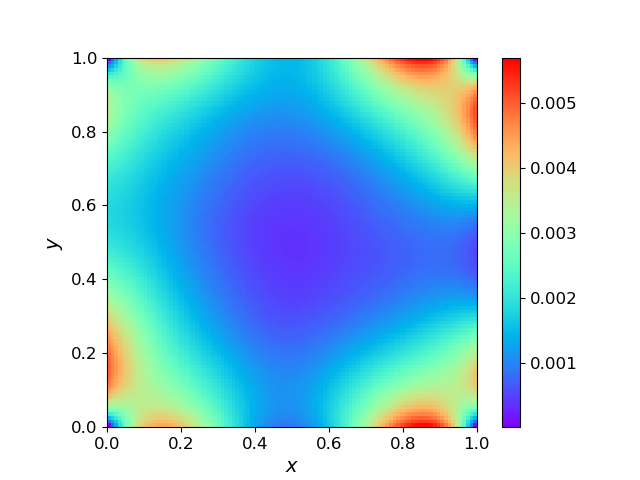} \\
(a) $\rho^{\text{Ref}}$, $t=0.1$ & (b) $\rho^{nn}_\theta$, $t=0.1$ & (c) $|\rho^{\text{Ref}}-\rho^{nn}_\theta|$, $t=0.1$ \\
\includegraphics[width=0.32\textwidth]{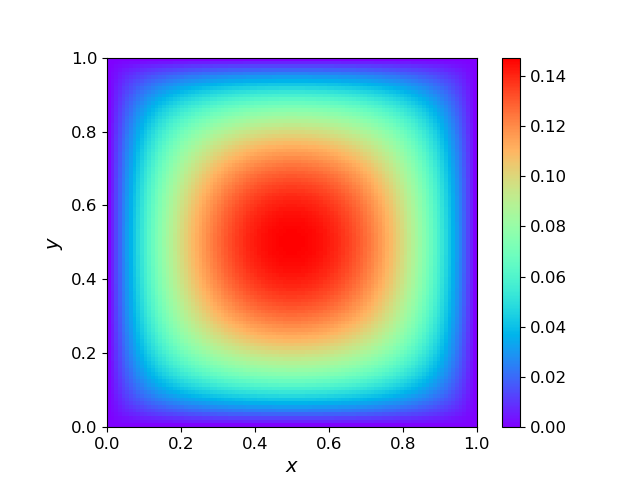} &
\includegraphics[width=0.32\textwidth]{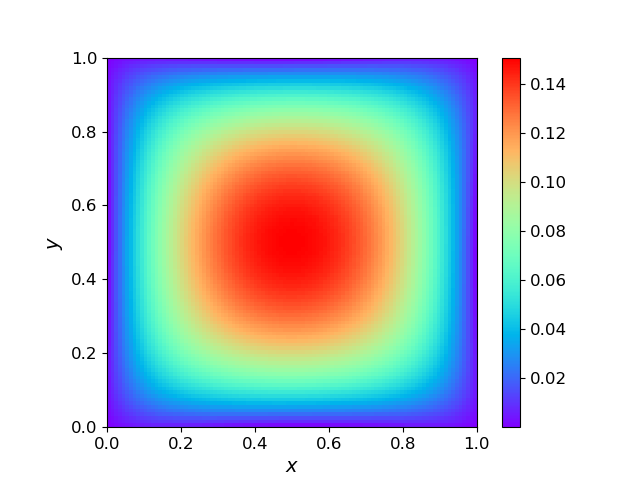} &
\includegraphics[width=0.32\textwidth]{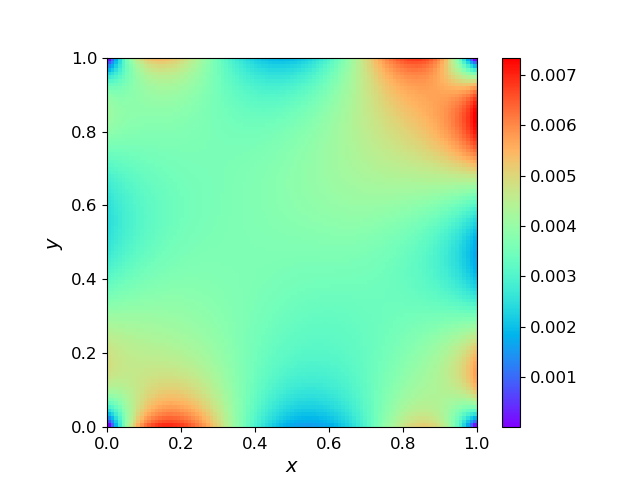} \\
(d) $\rho^{\text{Ref}}$, $t=0.8$ & (e) $\rho^{nn}_\theta$, $t=0.8$ & (f) $|\rho^{\text{Ref}}-\rho^{nn}_\theta|$, $t=0.8$ \\
\end{tabular}
\caption{The 2D LRTEs in the diffusion regime ($\epsilon=10^{-8}$).}
\label{figure4.2.2}
\end{figure}

We use the same network structure and parameters as the previous example except that $\beta_1=10^{-5}, \beta_2=10^{-16}$.
The reference solution is obtained by the finite difference discretization of the diffusion limit equation. We show the reference solution and the prediction of MA-APNNs in Figure~\ref{figure4.2.2}. The $L^2$ relative errors at times $t=0.1, 0.8$ are $3.98e-02$ and $4.79e-02$, respectively.

The above two examples show the ability of MA-APNNs to solve the 2D LRTEs in both the kinetic and diffusion regimes.
\subsection{\textbf{Uncertainty quantification (UQ) problems}}\label{ex:10}
Let $\{z_i\}_{i=1}^{d}$ be the random variables with identical independent uniform distribution in $[-1,1]$, i.e., $z_i \sim U([-1,1])$.
We apply the MA-APNNs to solve the uncertainty quantification problem \eqref{eq:LRTEs-uq} in the spatial, temporal and angular domain $D\times \tau \times S^2 =(0,1) \times (0,1) \times [-1,1]$, with $\alpha=0$, and the zero initial conditions,
\[
\begin{aligned}
f(0,x,\mu,\bm{z})=0.
\end{aligned}
\]
The model is a spatial 1D LRTEs with the high-dimensional random input $\bm{z} \sim U([-1,1]^{d})$.
The simulation is conducted for two cases with various scattering coefficients $\sigma(\bm{z})$.

\begin{figure}[!htbp]
	\centering
	{
		\begin{minipage}{2.3in}
			\centering
			\includegraphics[width=2.3in]{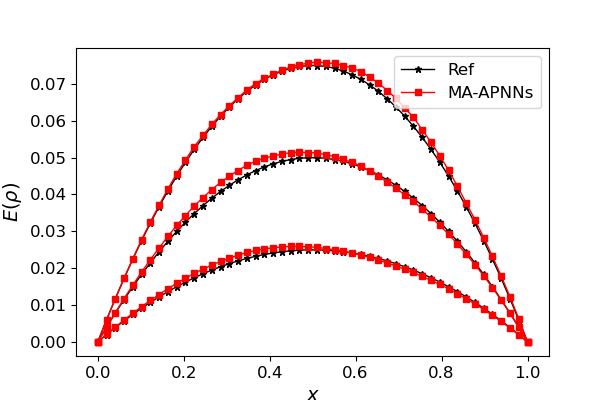}
		\end{minipage}
	}
	\caption{The expectation $\bm{E}(\rho)$ at $t=0.2, 0.4, 0.6$.}
	\label{figure4.3.2}
\end{figure}

\textbf{Problem \uppercase\expandafter{\romannumeral1}.} We consider the case in the kinetic regime $(\epsilon = 1)$ with the cosine scattering coefficient:
\[
\sigma(\bm{z})=1+0.1\sum \limits_{i=1}^{10}\cos(\pi z_i), \quad (d=10).
\]
The source term is set as
\[
G=\frac{x(1-x)}{22}\Big(\mu + 11 + \sum \limits_{i=1}^{10}z_i \Big) + \frac{\mu t (1-2x) }{22 \epsilon} \Big(\mu + 11 + \sum \limits_{i=1}^{10}z_i \Big) + \frac{1}{\epsilon^2}\sigma(\bm{z})tx(1-x)\mu.
\]
We impose the boundary conditions
\[
\begin{aligned}
f(t,0,\mu, \bm{z})= f(t,1,\mu,\bm{z})=0, \quad \mu \in [-1,1].
\end{aligned}
\]
The exact solution is given by
\[
f(t,x,\mu, \bm{z})=\frac{tx(1-x)}{22}\Big( \mu + 11 + \sum \limits_{i=1}^{10}z_i \Big), \quad  \rho(t,x,\bm{z})=\frac{tx(1-x)}{22}\Big(11 + \sum \limits_{i=1}^{10}z_i\Big).
\]
Note that the expectation of $\rho$ with respect to the random variable $\bm{z}$ is
\[
\bm{E}(\rho)=\frac{1}{2}tx(1-x).
\]
We take $\tilde{f}_{\theta}^{nn}=[13,40,40,40,40,1]$, and the boundaries and initial conditions are enforced by hard constraints, i.e., $f_{\theta}^{nn}= t x (1-x)\tilde{f}_{\theta}^{nn}$.
And we set $N_{\text{int}}=5000$, $\lambda_{b}=\lambda_i=\lambda_{c}=0, \beta_1=10^{-5}, \beta_2=10^{-7}$.
We carry out $10^4$ times simulations for $(\bm{z_1},\bm{z_2},...,\bm{z_{10}})$ and compute the expectation of $\rho$ at $t=0.2, 0.4, 0.6$.
The results are plotted in Figure~\ref{figure4.3.2}, which agrees well with the exact $\bm{E}(\rho)$.
Table~\ref{Table4.3.2} shows the $L^2$ errors of $\bm{E}(\rho)$.

\begin{table}[!htbp]
	\centering
	\caption{High dimensional transport regime. The errors of MA-APNNs.}
	\label{Table4.3.2}
	\begin{tabular}{llll}
		\hline\noalign{\smallskip}
		$L_{\text{error}}^2$ &$t=0.2$ & $t=0.4$ & $t=0.6$ \\
		\noalign{\smallskip}\hline\noalign{\smallskip}
		MA-APNNs &4.49e-02    &3.63e-02  &2.12e-02   \\
		\noalign{\smallskip}\hline
	\end{tabular}
\end{table}

\begin{figure}[!htbp]
	\centering
	{
		\begin{minipage}{2.3in}
			\centering
			\includegraphics[width=2.3in]{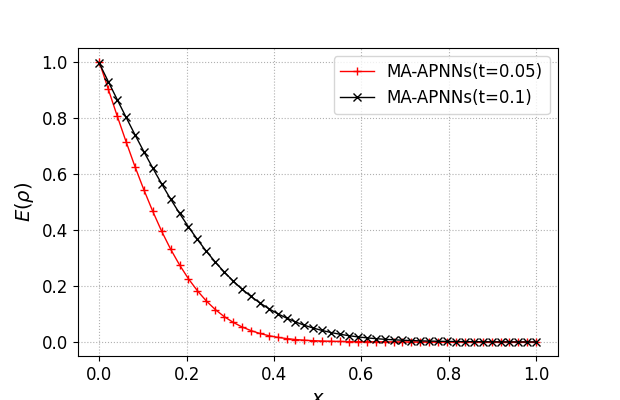}
		\end{minipage}
	}
	\caption{Plot of the density $\rho$ by taking expectation for $\bm{z}$ at $t=0.05, \ 0.1$. MA-APNNs with $\epsilon=10^{-5}$.}
	\label{figure4.3.22}
\end{figure}

\textbf{Problem \uppercase\expandafter{\romannumeral2}.}
We consider the case in the diffusion regime ($\epsilon = 10^{-5}$) with the sine scattering function\cite{jin2023asymptotic2}:
\[
\sigma(\bm{z})=1+0.1\prod \limits_{i=1}^{20}\sin(\pi \bm{z}_i), \quad (d=20).
\]
We impose the inflow boundary conditions and the zero source term,
\[
\begin{aligned}
f(t,0,\mu>0, \bm{z})=1, \quad f(t,1,\mu<0,\bm{z})=0, \quad G=0.
\end{aligned}
\]
We choose $f_{\theta}^{nn}=[23,40,40,40,40,1]$, $\left( N_{\text{int}},N_{\text{sb}}, N_{\text{tb}} \right)=(2048, 768\times 2, 1536)$ and set $(\beta_1,\beta_2,\lambda_{b},\lambda_i,\lambda_{c})= (10^{-5}, 10^{-16}, 1,1,0)$ in our simulation. We evaluate the expectation of $\rho$ at times $t=0.05, 0.1$. The results are shown in Figure~\ref{figure4.3.22}, which is very similar to the solutions obtained by \cite{jin2023asymptotic2}.
The above two examples indicate the applicability of MA-APNNs to the 1D LRTEs with high-dimensional random input in both kinetic and diffusion regimes.

\section{Conclusion}\label{sec5}
In this work, we derive the macroscopic auxiliary equation of the linear radiative transfer equations (LRTEs) and propose the macroscopic auxiliary asymptotic-preserving neural networks (MA-APNNs) to tackle the multiscale characteristics and high dimensionality in the numerical computation of LRTEs. To validate the advantages of MA-APNNs compared with the vanilla PINNs and APNNs in \cite{jin2023asymptotic1}, we conduct the simulation of the 1D LRTEs in kinetic, diffusive, and intermediate regimes.
We also present the numerical examples of the 2D LRTEs and 1D LRTEs with high-dimensional random variables, in both kinetic and diffusive regimes, showing the applicability of MA-APNNs to deal with the multiscale nature and high dimensionality.
Our future study will focus on developing MA-APNNs to solve the Vlasov-Poisson-Fokker-Planck system in the high-field regime.

\section*{Acknowledgements}
	This work was supported by the National Key R\&D Program (2020YFA0712200) for S. Jiang, and by Beijing Natural Science Foundation(Grant No. Z230003) for W. Sun, by NSFC (No. 12071060 \& No. 62231016) for L. Xu and by NSFC General Projects (No. 12171071) for G. Zhou.

	\bibliography{MAAPNN}
\end{document}